\documentclass{article}
\usepackage{amsmath}
\usepackage{amsfonts}
\usepackage{mathrsfs}
\usepackage[dvips]{graphics}
\input amssymb.sty
\begin{document}
\newcommand{\id}{\indent}
\newcommand{\nid}{\noindent}
\newcommand{\nn}{\nonumber}
\renewcommand{\ss}{\displaystyle}

\newenvironment{ea}{\begin{eqnarray}}{\end{eqnarray}}
\newenvironment{ea*}{\begin{eqnarray*}}{\end{eqnarray*}}
\newenvironment{ar}{\begin{array}}{\end{array}}
\newenvironment{eq}{\begin{equation}}{\end{equation}}

\baselineskip18pt \addtocounter{page}{0}
\title{\Large\bf Multi-almost periodicity and invariant basins of general neural networks under almost periodic stimuli
\thanks{This work was supported by the Foundation of Education of Fujian
Province, China (JA07142), the Foundation for Young Professors of
Jimei University, the Scientific Research Foundation of Jimei
University, the Foundation for Talented Youth with Innovation in
Science and Technology of Fujian Province (2009J05009).\newline\dag
E-mail: hzk974226@jmu.edu.cn}}
\author{Zhenkun Huang$^{*,\dag}$ \\
{\small {\sl School of Sciences, Jimei University, Xiamen 361021,
China } }}
\date{}
\maketitle \small \baselineskip13pt \begin{center}\textmd{SUMMARY}
\end{center} In this paper, we investigate convergence
dynamics of $2^N$ almost periodic encoded patterns of general neural
networks (GNNs) subjected to external almost periodic stimuli,
including almost periodic delays. Invariant regions are established
for the existence of $2^N$  almost periodic encoded patterns under
two classes of activation functions. By employing the property of
$\mathscr{M}$-cone and inequality technique, attracting basins are
estimated and some criteria are derived for the networks to converge
exponentially toward $2^N$ almost periodic encoded patterns. The
obtained results are new, they extend and generalize the
corresponding results existing in previous literature.

\vskip 0.15in\noindent {\scriptsize KEY WORDS:}\ \ neural networks;
almost periodic; encoded patterns; attracting basins;  exponential
stability

\vskip 0.5in \normalsize \baselineskip16pt
\centerline{1.~~\textrm{INTRODUCTION}}

\vskip 0.1in\noindent In recent years, the dynamical behaviors of
neural networks with delays have been widely investigated. Many
important results on the existence and uniqueness of equilibrium
point, global asymptotic (exponential) stability have been
established and successfully applied to signal and image processing
system, associative memories, pattern classification and so on. For
corresponding results, we can refer to [1-9,17-18,34-36]. In the
applications of neural networks to associative memory storage or
pattern recognitions, the existence of multiple equilibria or
multiple periodic orbits is an important feature [1-3,6-9]. It is
worth noting that convergence analysis and coexistence of multiple
equilibria or multiple periodic solutions have been investigated in
[6-8] and these equilibria or periodic solutions are usually called
\textbf{encoded patterns} [6-8,10].

As we know well, the nonautonomous phenomenon involved in periodic
or almost periodic environment often occurs in many realistic
systems [11-14,19]. Hence, in many applications, the property of
periodic or almost periodic oscillatory solutions of neural networks
is of great interest. Meanwhile, there often exist delays in
artificial neural networks due to the finite switching speed of
amplifiers and faults in the electrical circuit. They slow down the
transmission rate and lead to some degree of instability. Therefore,
complex dynamic behaviors of neural networks under periodic or
almost periodic stimuli and delayed effects have been studied so far
[10,15-16,20-31,37,39,41].

However, to the best of the authors' knowledge, few papers  deal
with general neural networks with both almost periodic coefficients
and almost periodic delays. Furthermore, most of the results
reported in the literature focus on the stability of unique almost
periodic solution of neural networks. We can refer to [10,20-23,39]
and the references cited therein. In this paper, similarly as
[26,31], we consider the following nonautonomous general neural
networks with transmission delays
\begin{ea}
\frac{du^i(t)}{dt}=\displaystyle-c_i(t)u^i(t)+\sum\limits_{l=1}^M
\sum\limits_{j=1}^Na_{ijl}(t)g_j\Big(\sigma_ju^j\big(t-\kappa_{ijl}(t)\big)\Big)+J_i(t),
\end{ea}
where $i$, $j$$\in$$\mathscr{I}$:=$\{1,2,\cdots,N\}$,
$l$$\in$$\mathscr{L}$:=$\{1,2,\cdots,M\}$; The main purpose of this
paper is to study complex convergence dynamics of GNNs (1) in
encoding external stimuli that vary almost periodically with time
and recalling the encoded patterns associated with almost periodic
delays. That is, we investigate  exponential stability of $2^N$
almost periodic encoded patterns (almost periodic solutions) of GNNs
(1). The criteria we established are completely different from most
of the existing results in [10,20-31]. Particularly, when GNNs (1)
degenerates into the autonomous system, our results extend and
generalize the related results existing in [6,8].

The rest of this paper is organized as follows. In Section 2, we
shall make some preparations by giving some definitions and lemmas.
Meanwhile, we establish $2^N$ positively invariant basins for
general neural networks under almost periodic stimuli. In Section 3,
by using the property of $\mathscr{M}$-cone and inequality
technique, attracting basins are determined and some new criteria
for exponential stability of $2^N$ almost periodic encoded patterns
are obtained. In Section 4, we shall make some generalizations by
considering activation functions with saturations and apply our
obtained results to some special neural networks systems. It is
shown that our results are general and improve the previously known
results. Finally, numerical examples are presented to illustrate our
results.

\vskip 0.3in \normalsize \baselineskip16pt
\centerline{2.~~\textrm{PRELIMINARIES AND $2^N$ INVARIANT BASINS FOR
GNNS}} \vskip 0.1in\noindent In GNNs (1), the integer N corresponds
to the number of units in neural networks and M corresponds to the
number of neural axons; that is, signals that emit from the $i$th
unit have M pathways to the $j$th unit; $u^i(t)$ corresponds to the
membrane potential of the $i$th neuron at time $t$; the dissipation
coefficient $c_i(t)>0$ represents the rate with which the $i$th
neuron resets its potential when isolated from other neurons and
inputs; $a_{ijl}(t)$ denotes synaptic connection weight of the $j$th
neuron on the $i$th neuron at time $t-\kappa_{ijl}(t)$; $J_i(t)$ is
an input to the $i$th neuron at time $t$ from outside the networks;
$g_{j}(\cdot)$ denotes activation function; $\sigma_j$ denotes the
amplifier gain; $\kappa_{ijl}(t)$ is the transmission delay of the
$i$th unit along the $l$ axon of the $j$th unit at time t, it is a
nonnegative bounded function with
$0<\kappa_{ijl}(t)\leq\kappa_{ijl}<\kappa:=\max\{\kappa_{ijl}|
i,j\in\mathscr{I},l\in\mathscr{L}\}$.

As usual, we denote by ${\mathcal{C}}([-\kappa,0],\mathscr{R}^N)$
the Banach space of all real-valued continuous mappings from
$[-\kappa,0]$ to $\mathscr{R}^N$ equipped with supremum norm defined
by \begin{ea*} \|\phi\|_\kappa=\max\limits_{1\leq i\leq
N}\|\phi^i\|_\kappa,\,\,\,\,
\mbox{where}\,\,\,\,\|\phi^i\|_\kappa=\sup\limits_{-\kappa\leq
t\leq0}|\phi^i(t)| \end{ea*} and
$\phi=(\phi^1,\phi^2,\cdots,\phi^N)^T\in{\mathcal{C}}([-\kappa,0],\mathscr{R}^N)$.
Let $l>0$. For any $u\in{\mathcal{C}}([-\kappa,l],\mathscr{R}^N)$
and $t\in[0,l]$, we define $ u_t(s)=u(t+s)$, $s\in[-\kappa,0]$. Then
we have $u_t(\cdot)\in{\mathcal{C}}([-\kappa,0],\mathscr{R}^N)$. For
any given $\phi\in{\mathcal{C}}([-\kappa,0],\mathscr{R}^N)$, we
denote by $u(t;\phi)$ the solution of GNNs (1) with $u_0(s)=\phi(s)$
for all $s\in[-\kappa,0]$.

\vskip 0.1in\noindent\textsl{Definition 1} (see [12-13])\\ A
continuous function $f(t):\mathscr{R}\rightarrow\mathscr{R}$ is
 called an almost periodic function if for any $\epsilon>0$,
$$\mathscr{T}(f,\epsilon)=\Big\{\mathcal{T}\in \mathscr{R}\Big|
|f(t+\mathcal{T})-f(t)|<\epsilon\ \ \mbox{for all}\ \ t\in
\mathscr{R}\Big\}$$ is a relatively dense set in $\mathscr{R}$. That
is, there exists a positive constant $l(\epsilon)$ such that any
interval with length $l(\epsilon)$ contains at least one point of
$\mathscr{T}(f,\epsilon)$. The number $\mathcal{T}$ is called
$\epsilon$-almost period of $f(t)$.

\vskip 0.1inLet $(\mathscr{AP},\|\cdot\|)$ be the space of all
real-valued almost periodic functions defined on $\mathscr{R}$ with
supremum norm defined by $\|f\|=\sup\limits_{t\in\mathscr{R}}|f(t)|$
for any $f(t)\in\mathscr{AP}$. It is easy for us to have the
following basic properties (see [12-13]):

\vskip 0.1in $(\blacktriangle_1)$ Let $f(t)\in\mathscr{AP}$. Then
$f(t)$ is bounded and uniformly continuous for all
$t\in\mathscr{R}$.

\vskip 0.1in $(\blacktriangle_2)$ Let $f_i(t)\in\mathscr{AP}$,
$i=1,2,\cdots,N$. Then for any $\epsilon>0$, there exists a constant
$l(\epsilon)>0$ such that any interval with length $l(\epsilon)$
contains at least one common $\epsilon$-almost period of $f_i(t)$
for all $i=1,2,\cdots,N$. That is,
$\bigcap\limits_{i=1}^N\mathscr{T}(f_i,\epsilon)$ is relatively
dense in $\mathscr{R}$.

\vskip 0.1in $(\blacktriangle_3)$ Let $f(t),g(t)\in\mathscr{AP}$.
Then $f(t-g(t))\in\mathscr{AP}$.

\vskip 0.1in\noindent \textsl{Remark 1}\\ Assume that
$f(t),g(t)\in\mathscr{AP}$. For any $\epsilon>0$, by using
$(\blacktriangle_1)$, there exists $\delta=\delta(\epsilon/2)>0$
such that $|f(t_1)-f(t_2)|<\epsilon/2$ if $|t_1-t_2|<\delta$. We
only take $\widehat{\epsilon}=\min\{\epsilon/2,\delta\}$. It is easy
for us to check that
$|f(t+\mathcal{T}-g(t+\mathcal{T}))-f(t-g(t))|<\epsilon$ for any
$\mathcal{T}\in\mathscr{T}\big(f,\widehat{\epsilon}\,\big)\bigcap\mathscr{T}\big(g,\widehat{\epsilon}\,\big)$
and all $t\in\mathscr{R}$. From $(\blacktriangle_2)$ and Definition
1, we get $f(t-g(t))\in\mathscr{AP}$. That is, the property
$(\blacktriangle_3)$ holds.

\vskip 0.1in \noindent Throughout this paper, we always assume the
following assumptions hold: \vskip 0.1in$\bullet$\,$(H_1)$\
$c_i(t)>0$, $a_{ijl}(t)$, $\kappa_{ijl}(t)>0$ and $J_i(t)$ are
almost periodic functions defined on $\mathscr{R}$, $\sigma_j>0$ and
$a_{iil}(t)>0$, where $i,j$$\in\mathscr{I}$ and $l\in\mathscr{L}$.

\vskip 0.1in\noindent Unless otherwise stated, we always use
$i,j=1,2,\cdots,N$, $l=1,2,\cdots,M$. The activation functions of
first class we considered satisfy the following basic assumption:
\begin{ea*}
\mbox{Class}\,\mathcal{A}\,:\,\,\, g_i\in\mathcal{C}^2,\,\,\,
\left\{\begin{aligned}
|g_i(x)|&\leq B_i,\,\,g_i(0)=0\,\,\, \mbox{and} \\
\dot{g}_i(x)&>0,\,\,\,x\ddot{g}_i(x)<0, \,\,\,\mbox{where} \,\,\,
x\in\mathscr{R}.
\end{aligned} \right.
\end{ea*}

\noindent\textsl{Lemma 1}\\Assume that $(H_1)$ holds. For any
$\phi=(\phi^1,\phi^2,\cdots,\phi^N)^T\in{\mathcal{C}}([-\kappa,0],\mathscr{R}^N)$,
$$\|\phi^i\|_\kappa\leq\Big(\sum\limits_{l=1}^M\sum\limits_{j=1}^N\sup\limits_{t\in\mathscr{R}}|a_{ijl}(t)|B_j
+\sup\limits_{t\in\mathscr{R}}|J_i(t)|\Big)
\Big/\inf\limits_{t\in\mathscr{R}}c_i(t)$$ implies that
$$\|u_t^i(\cdot;\phi)\|_\kappa\leq\Big(\sum\limits_{l=1}^M\sum\limits_{j=1}^N\sup\limits_{t\in\mathscr{R}}|a_{ijl}(t)|B_j
+\sup\limits_{t\in\mathscr{R}}|J_i(t)|\Big)
\Big/\inf\limits_{t\in\mathscr{R}}c_i(t)$$ for all $t\geq 0$, where
$u(t;\phi)$ is the solution of GNNs (1) with $u_0(s)=\phi(s)$ for
$s\in[-\kappa,0]$.

\vskip 0.1in\noindent \textsl{Proof}\\The proof is trivial, we omit
it here.
\indent\indent\indent\indent\indent\indent\indent\indent\indent\indent\indent\indent\indent\indent\indent\indent\indent\indent\indent$\square$

\vskip 0.1in We introduce the following auxiliary functions
$$F_i(z)=-\sup\limits_{t\in\mathscr{R}}c_i(t)z+\inf\limits_{t\in\mathscr{R}}\sum\limits_{l=1}^Ma_{iil}(t)g_i(\sigma_iz).$$

\noindent\textsl{Lemma 2}\\Suppose that the following
assumption holds:\\
\indent$\bullet(H_1^\mathcal{A})\,\,\inf\limits_{\zeta\in\mathscr{R}}\dot{g}_i(\zeta)<\ss
\sup\limits_{t\in\mathscr{R}}c_i(t)\Big/\sigma_i\inf\limits_{t\in\mathscr{R}}\sum\limits_{l=1}^Ma_{iil}(t)<
\sup\limits_{\zeta\in\mathscr{R}}\dot{g}_i(\zeta).$\\ Then there
exist two points $z_{i1}$ and $z_{i2}$ with $z_{i1}<0<z_{i2}$ such
that $\dot{F}_i(z_{ik})=0$ and
$\dot{F}_i(z){\cdot}sgn\Big\{\displaystyle\frac{z-z_{i1}}{z-z_{i2}}\Big\}<0$
$(z\neq z_{ik},k=1,2)$.

\vskip 0.1in\noindent \textsl{Proof}\\We have $\dot{F}_i(z)=0$ if
and only if
$\ss\dot{g}_i(\sigma_iz)=\sup\limits_{t\in\mathscr{R}}c_i(t)\Big/\sigma_i\inf\limits_{t\in\mathscr{R}}\sum\limits_{l=1}^Ma_{iil}(t).$
For activation functions of class $\mathcal{A}$, we know that the
graph of positive function $\dot{g}_i(z)$ concaves down and has its
maximal value at zero. By the continuity of $\dot{g}_i(z)$ and
$(H_1^\mathcal{A})$, there exist two points $z_{i1}$ and $z_{i2}$
with $z_{i1}<0<z_{i2}$ such that
$\dot{g}_i(\sigma_iz_{ik})=\sup\limits_{t\in\mathscr{R}}c_i(t)\Big/\sigma_i\inf\limits_{t\in\mathscr{R}}\sum\limits_{l=1}^Ma_{iil}(t)$;
that is, $\dot{F}_i(z_{ik})=0$ $(k=1,2)$. Since $\dot{g}_i(z)$ is
increasing on $(-\infty,z_{i1}]$ and is decreasing on
$[z_{i2},+\infty)$, we get that
$$\Big(-\sup\limits_{t\in\mathscr{R}}c_i(t)+\sigma_i\dot{g}_i(\sigma_iz)\inf\limits_{t\in\mathscr{R}}\sum\limits_{l=1}^Ma_{iil}(t)
\Big){\cdot}sgn\Big\{\displaystyle\frac{z-z_{i1}}{z-z_{i2}}\Big\}<0;$$
that is,
$\dot{F}_i(z){\cdot}sgn\Big\{\displaystyle\frac{z-z_{i1}}{z-z_{i2}}\Big\}<0$
$(z\neq z_{ik})$. The proof is complete.
\indent\indent\indent\indent\indent\indent\indent$\square$

\vskip 0.1in For the existence of $2^N$ positively invariant basins
of GNNs (1), we consider the following assumption for activation functions of class $\mathcal{A}$:\\
\indent$\bullet(H_2^\mathcal{A})\,\,(-1)^k\cdot\Big\{F_i(z_{ik})
+J_i(t)\Big\}>\ss\sum\limits_{l=1}^M\sum\limits_{j\neq{i}}\sup\limits_{t\in\mathscr{R}}|a_{ijl}(t)|B_j$\\
for all $t\in\mathscr{R}$, where $k=1,2$. \vskip 0.1in \indent
Take $k=1$ in  $(H_2^\mathcal{A})$, it is easy for us to get that
\begin{ea}F_i(z_{i1})+\sum\limits_{l=1}^M\sum\limits_{j\neq{i}}\sup\limits_{t\in\mathscr{R}}|a_{ijl}(t)|B_j
+\sup\limits_{t\in\mathscr{R}}J_i(t)<0.
\end{ea} Noting that
$F_i(z)\rightarrow+\infty$ as $z\rightarrow-\infty$, we know that
there exists a $\widehat{z}_{i1}$ with $\widehat{z}_{i1}<z_{i1}<0$
such that
\begin{ea}
F_i(\widehat{z}_{i1})+\sum\limits_{l=1}^M\sum\limits_{j\neq{i}}\sup\limits_{t\in\mathscr{R}}|a_{ijl}(t)|B_j
+\sup\limits_{t\in\mathscr{R}}J_i(t)=0.
\end{ea} Take $k=2$ in
$(H_2^\mathcal{A})$, by the similar argument, we derive that there
exists a $\widetilde{z}_{i2}$ with $0<z_{i2}<\widetilde{z}_{i2}$
such that
\begin{ea}
F_i(\widetilde{z}_{i2})-\sum\limits_{l=1}^M\sum\limits_{j\neq{i}}\sup\limits_{t\in\mathscr{R}}|a_{ijl}(t)|B_j
+\inf\limits_{t\in\mathscr{R}}J_i(t)=0.
\end{ea}
Next, we let
\begin{ea*}
\left\{\begin{ar}{l}\ss\alpha_{i1}:=-\Big(\sum\limits_{l=1}^M\sum\limits_{j=1}^N\sup\limits_{t\in\mathscr{R}}|a_{ijl}(t)|B_j
+\sup\limits_{t\in\mathscr{R}}|J_i(t)|\Big)\Big/\inf\limits_{t\in\mathscr{R}}c_i(t),\,\,\,\beta_{i1}:=\widehat{z}_{i1},
_{\left.\begin{ar}{l}\\\end{ar}\right.}\\
\alpha_{i2}:=\widetilde{z}_{i2},\,\,\,\,\beta_{i2}:=\Big(\sum\limits_{l=1}^M\sum\limits_{j=1}^N\sup\limits_{t\in\mathscr{R}}|a_{ijl}(t)|B_j
+\sup\limits_{t\in\mathscr{R}}|J_i(t)|\Big)\Big/\inf\limits_{t\in\mathscr{R}}c_i(t).
\end{ar}\right.
\end{ea*}
It is easy for us to check that
$\alpha_{i1}<\beta_{i1}<0<\alpha_{i2}<\beta_{i2}$. Then we define
$2N$ subsets of ${\mathcal{C}}([-\kappa,0],\mathscr{R})$ as follows:
\begin{ea*}
&&\mathscr{K}_{i1}:=\Big\{\psi\in{\mathcal{C}}([-\kappa,0],\mathscr{R})\,
|\, \psi(s)\leq\beta_{i1} \ \mbox{for}\ \mbox{all}\
s\in[-\kappa,0]\Big\},\\
&&\mathscr{K}_{i2}:=\Big\{\psi\in{\mathcal{C}}([-\kappa,0],\mathscr{R})\,
|\, \psi(s)\geq\alpha_{i2} \ \mbox{for}\ \mbox{all}\
s\in[-\kappa,0]\Big\}.
\end{ea*}
Hence we have $2^N$ subsets $
\mathscr{K}^\Sigma:=\underbrace{\mathscr{K}_{1\varsigma_1}\times\mathscr{K}_{2\varsigma_2}
\cdots\times\mathscr{K}_{N\varsigma_N}}_N $ of
${\mathcal{C}}([-\kappa,0],\mathscr{R}^N)$, where
$\Sigma=(\varsigma_1,\varsigma_2,\cdots,\varsigma_N)$ with
$\varsigma_i=1$ or $2$, $i\in\mathscr{I}$. In what follows, we
should prove that these $ \mathscr{K}^\Sigma$ are $2^N$ positively
invariant basins of GNNs (1).

\vskip 0.1in\noindent\textsl{Theorem 1}\\ Under the assumptions
$(H_1)$ and $(H_1^\mathcal{A}$)-(H$_2^\mathcal{A}$), each
$\mathscr{K}^\Sigma$ is a positively invariant basin with respect to
the solution flow generated by GNNs (1).

\vskip 0.1in\noindent \textsl{Proof}\\For any initial condition
$\phi\in\mathscr{K}^\Sigma$, we should prove that
$u_t(\cdot;\phi)\in\mathscr{K}^\Sigma$ for all $t\geq0$. For each
$i\in\mathscr{I}$, we only consider the case $\varsigma_i=2$, i.e.,
$\phi^i(s)\geq\alpha_{i2}$ for all $s\in[-\kappa,0]$. We assert
that, for any sufficiently small $\epsilon>0$ $(\epsilon\ll
\alpha_{i2}-z_{i2})$, the solution
$u^i(t;\phi)\geq\alpha_{i2}-\epsilon$ holds for all $t\geq0$. If
this is not true, there exists some $t^*>0$ such that
$u^i(t^*)=\alpha_{i2}-\epsilon$, $\dot{u}^i(t^*)\leq0$ and
$u^i(t)>\alpha_{i2}-\epsilon$ for $t\in[-\kappa,t^*]$. Due to
$(H_1^\mathcal{A}$),
$\inf\limits_{t\in\mathscr{R}}\sum\limits_{l=1}^Ma_{iil}(t)>0$ and
the monotonicity of $g_{i}$, we derive from GNNs (1) that
\begin{ea}
\frac{du^i(t^*)}{dt}\hskip-0.1in&=&\hskip-0.1in\displaystyle-c_i(t^*)u^i(t^*)
+\sum\limits_{l=1}^M\sum\limits_{j=1}^Na_{ijl}(t^*)g_j\Big(\sigma_ju^j\big(t^*-\kappa_{ijl}(t^*)\big)\Big)
+J_i(t^*)\nn\\
\hskip-0.1in&\geq&\hskip-0.1in-\sup\limits_{t\in\mathscr{R}}c_i(t)(\alpha_{i2}-\epsilon)+
\inf\limits_{t\in\mathscr{R}}\sum\limits_{l=1}^Ma_{iil}(t)g_i\big(\sigma_i(\alpha_{i2}-\epsilon)\big)\nn\\
\hskip-0.1in&-&\hskip-0.1in\sum\limits_{l=1}^M\sum\limits_{j\neq{i}}\sup\limits_{t\in\mathscr{R}}|a_{ijl}(t)|B_j
+\inf\limits_{t\in\mathscr{R}}J_i(t)\nn\\
\hskip-0.1in&\geq&\hskip-0.1inF_i(\widetilde{z}_{i2}-\epsilon)-
\sum\limits_{l=1}^M\sum\limits_{j\neq{i}}\sup\limits_{t\in\mathscr{R}}|a_{ijl}(t)|B_j
+\inf\limits_{t\in\mathscr{R}}J_i(t).
\end{ea}\nid
From Lemma 2, we know that $F_i(z)$ is strictly decreasing on
$(z_{i2},+\infty)$. By using (4) and (5), we get
$\frac{du^i(t^*)}{dt}>0$ which leads to a contradiction. Since the
choice of $\epsilon$ is arbitrary, if $\phi^i(s)\geq\alpha_{i2}$ for
all $s\in[-\kappa,0]$, we have $u^i(t;\phi)\geq\alpha_{i2}$ for all
$t\geq0$. When $\varsigma_i=1$, similar argument can be performed to
show that if $\phi^i(s)\leq\beta_{i1}$ for all $s\in[-\kappa,0]$, we
have $u^i(t;\phi)\leq\beta_{i1}$ for all $t\geq0$. Hence, for any
initial condition $\phi\in\mathscr{K}^\Sigma$, we have that
$u_t(\cdot,\phi)\in\mathscr{K}^\Sigma$ for all $t\geq0$. That is,
each $\mathscr{K}^\Sigma$ is a positively invariant basin with
respect to the solution flow generated by GNNs (1).  The proof is
complete. \indent\indent\indent\indent\indent\indent\indent$\square$

\vskip 0.1in \indent For convenience of discussing invariant regions
for the existence of $2^N$ almost periodic encoded patterns of GNNs
(1) in next section, we define $2N$ subsets
${\mathscr{B}}_{ik}\subset{\mathscr{AP}}$ ($i\in\mathscr{I}$,
$k=1,2$) which satisfy the following two basic properties:

\vskip 0.1in$(\bigtriangleup_1)$ For any
$\phi(t)\in{\mathscr{B}}_{ik}$,
$\alpha_{ik}\leq\phi(t)\leq\beta_{ik}$ for all $t\in\mathscr{R}$.

\vskip 0.1in$(\bigtriangleup_2)$ For any $\epsilon>0$,
${\mathscr{T}}({\mathscr{B}}_{ik},\epsilon):=\bigcap
\limits_{\phi\in{\mathscr{B}}_{ik}}{\mathscr{T}}(\phi,\epsilon)$ is
relatively dense.

\vskip 0.1in\noindent Let
\begin{ea*}
\mathscr{B}^\Sigma:=\underbrace{\mathscr{B}_{1\varsigma_1}\times\mathscr{B}_{2\varsigma_2}
\cdots\times\mathscr{B}_{N\varsigma_N}}_N\,\,,
\end{ea*}
where $\Sigma=(\varsigma_1,\varsigma_2,\cdots,\varsigma_N)$ with
$\varsigma_i=1$ or $2$, $i\in\mathscr{I}$. Then ${\mathscr{B}}_{ik}$
($i\in\mathscr{I}$, $k=1,2$) are not only convex subsets of Banach
space ${\mathscr{AP}}$, but also uniformly almost periodic families
(see [12]). The compactness of ${\mathscr{B}}_{ik}$ comes from the
following lemma.

\vskip 0.1in\noindent\textsl{Lemma 3} (see [12])\\ If
$\mathscr{B}\subset{\mathscr{AP}}$ is a uniformly almost periodic
family, then from every sequence in $\mathscr{B}$ one can extract a
subsequence which converges uniformly on $\mathscr{R}$.

\vskip 0.1in\noindent\textsl{Lemma 4} (see [26])\\For any $p>1$,
$x_k\geq0$, $y\geq0$, the following inequality holds:
$$y\prod\limits_{k=1}^mx_k^{p_k}\leq\frac{1}{p}\sum\limits_{k=1}^mp_kx_k^p+\frac{1}{p}y^p,$$
where $p_k>0$ ($k=1,2,\cdots,m$) are constants and
$\sum\limits_{k=1}^mp_k=p-1$.

\vskip 0.1in\noindent\textsl{Lemma 5} (see [32-34])\\ Let
$H=(h_{ij})_{N\times{N}}\in\mathscr{R}^{N\times{N}}$ with
$h_{ij}\leq0$ ($i\neq j$). Then the following conditions are
equivalent:\\
(1) All the leading principal minors of $H$ are positive;\\
(2)$H$ is quasi-dominant positive diagonal; that is, there exist
positive numbers $z_j$ ($j\in\mathscr{I}$) such that
$\sum\limits_{j=1}^Nz_jh_{ij}>0\,\,\, \mbox{or}
\,\,\,\sum\limits_{j=1}^Nz_jh_{ji}>0, \,\,\,i\in\mathscr{I}.$

\vskip 0.1inWe denote by $\mathscr{M}$ the set of all matrices which
satisfy one of the above properties. For any $H\in\mathscr{M}$, let
$\Omega_\mathscr{M}(H):=\Big\{Z=(z_1,z_2,\cdots,z_N)^T\in\mathscr{R}^{N}\Big|HZ>0\,\,\mbox{and}\,\,
z_i>0,i\in\mathscr{I}\Big\}.$ It is obvious that
$\Omega_\mathscr{M}(H)$ is a cone without the vertex in
$\mathscr{R}^{N}$. Given any
$\widetilde{H}=(\widetilde{h}_{ij}|\widetilde{h}_{ij}\leq 0,i\neq
j)_{N\times{N}}$. If $\widetilde{h}_{ij}\geq h_{ij}$
($i,j\in\mathscr{I}$) and $H\in\mathscr{M}$, then
$\widetilde{H}\in\mathscr{M}$.

\vskip 0.3in\normalsize \baselineskip16pt
\centerline{3.~~\textrm{ALMOST PERIODIC ENCODED PATTERNS FOR GNNS}}

\vskip 0.1in\noindent In this section, by using the properties of
almost periodicity and Schauder's fixed point theorem, we should
prove that each ${\mathscr{B}}^\Sigma$ is an invariant region and
there exist at least $2^N$ almost periodic encoded patterns of GNNs
(1) in these ${\mathscr{B}}^\Sigma$. Finally, attracting basins are
estimated and some  criteria are derived for the networks to
converge exponentially toward $2^N$ almost periodic encoded
patterns.

\vskip 0.1in\noindent\textsl{Theorem 2}\\Under the basic assumptions
$(H_1)$ and $(H_1^\mathcal{A})$-$(H_2^\mathcal{A})$, for each
$\Sigma$, there exists at least one almost periodic encoded pattern
of GNNs (1) in ${\mathscr{B}}^{\Sigma}$.

\vskip 0.1in\noindent \textsl{Proof}\\For each
$\Sigma=(\varsigma_1,\varsigma_2,\cdots,\varsigma_N)$, we define a
mapping
$\mathscr{F}^\Sigma=(\mathscr{F}_1^\Sigma,\mathscr{F}_2^\Sigma,\cdots,\mathscr{F}_N^\Sigma)$
by
\begin{ea}
({\mathscr{F}}_i^\Sigma\phi)(t)&\hskip-2mm=&\hskip-2mm\displaystyle\int\limits_0^\infty\Bigg[\sum\limits_{l=1}^M\sum\limits_{j=1}^N
a_{ijl}(t-s)g_j\Big(\sigma_j\phi^j(t-s-\kappa_{ijl}(t-s))\Big)\nn\\
&\hskip-2mm+&\hskip-2mmJ_i(t-s)\Bigg]\exp\Big(-\int\limits_0^sc_i(t-u)du\Big)ds,
\end{ea}where
$i\in\mathscr{I}$,
$\phi=(\phi^1,\phi^2,\cdots,\phi^N)\in{\mathscr{B}}^\Sigma$. From
$(H_1)$ and the boundedness of activation functions, it is easy for
us to check that each $\mathscr{F}^\Sigma$ is well defined. Next we
need three steps to complete our proof.

{\bf Step 1:}\,\,\,For each $i\in\mathscr{I}$, we should prove that
$\alpha_{i\varsigma_i}\leq({\mathscr{F}_i^\Sigma}\phi)(t)\leq\beta_{i\varsigma_i}$
for all $t\in\mathscr{R}$. Fix $i\in\mathscr{I}$. From $(H_1)$ and
(6), one obtains that
\begin{ea}
\Big|({\mathscr{F}}_i^\Sigma\phi)(t)\Big|\hskip-0.1in&\leq&\hskip-0.1in\displaystyle\Big(
\sum\limits_{l=1}^M\sum\limits_{j=1}^N\sup\limits_{t\in\mathscr{R}}|a_{ijl}(t)|B_j
+\sup\limits_{t\in\mathscr{R}}|J_i(t)|\Big)
\int\limits_0^\infty\exp\Big(-\int\limits_0^s\inf\limits_{t\in\mathscr{R}}c_i(t)du\Big)ds\nn
\end{ea}
\begin{ea}
\hskip-0.1in&\leq&\hskip-0.1in\Big(\sum\limits_{l=1}^M\sum\limits_{j=1}^N\sup\limits_{t\in\mathscr{R}}|a_{ijl}(t)|B_j
+\sup\limits_{t\in\mathscr{R}}|J_i(t)|\Big)\Big/\inf\limits_{t\in\mathscr{R}}c_i(t)=\beta_{i2}.
\end{ea}
If $\varsigma_i=2$ ($i\in\mathscr{I}$), then
$\phi^i(t)\geq\alpha_{i2}$ for all $t\in\mathscr{R}$. From (6) and
$(H_2^\mathcal{A})$, we get
\begin{ea}
({\mathscr{F}}_i^\Sigma\phi)(t)\hskip-0.1in&\geq&\hskip-0.1in\displaystyle\Big(\inf\limits_{t\in\mathscr{R}}
\sum\limits_{l=1}^Ma_{iil}(t)g_i(\sigma_i\alpha_{i2})
-\sum\limits_{l=1}^M\sum\limits_{j\neq{i}}\sup\limits_{t\in\mathscr{R}}|a_{ijl}(t)|B_j\Big)
\int\limits_0^\infty\exp\Big(-\int\limits_0^s\sup\limits_{t\in\mathscr{R}}c_i(t)du\Big)ds\nn\\
\hskip-0.1in&+&\hskip-0.1in\displaystyle\inf\limits_{t\in\mathscr{R}}J_i(t)
\int\limits_0^\infty\exp\Big(-\int\limits_0^s\sup\limits_{t\in\mathscr{R}}c_i(t)du\Big)ds\nn\\
\hskip-0.1in&\geq&\hskip-0.1in\Big(\inf\limits_{t\in\mathscr{R}}\sum\limits_{l=1}^Ma_{iil}(t)g_i(\sigma_i\alpha_{i2})
-\sum\limits_{l=1}^M\sum\limits_{j\neq{i}}\sup\limits_{t\in\mathscr{R}}|a_{ijl}(t)|B_j\Big)
\Big/\sup\limits_{t\in\mathscr{R}}c_i(t)\nn\\
\hskip-0.1in&+&\hskip-0.1in\displaystyle\inf\limits_{t\in\mathscr{R}}J_i(t)
\Big/\sup\limits_{t\in\mathscr{R}}c_i(t)=\alpha_{i2},
\end{ea}
for all $t\in\mathscr{R}$. By (7) and (8), we have
$\alpha_{i2}\leq({\mathscr{F}_i^\Sigma}\phi)(t)\leq\beta_{i2}$. From
similar argument, if $\varsigma_i=1$, we can prove that
$\alpha_{i1}\leq({\mathscr{F}}_i^\Sigma\phi)(t)\leq\beta_{i1}$ for
all $t\in\mathscr{R}$. Hence, we have
$\alpha_{i\varsigma_i}\leq({\mathscr{F}_i^\Sigma}\phi)(t)\leq\beta_{i\varsigma_i}$
for each $i\in\mathscr{I}$ and all $t\in\mathscr{R}$.

{\bf Step 2:}\,\,\,We should prove that
${\mathscr{F}^\Sigma}:{\mathscr{B}}^\Sigma\rightarrow{\mathscr{B}}^\Sigma$.
For any $\epsilon>0$, we let
\begin{ea}
\epsilon^*\hskip-0.1in&=\hskip-0.1in&\min\Big\{\,\,\frac{\epsilon}{5}\inf\limits_{t\in\mathscr{R}}c_i(t)
\Big/\max\Big[\sum\limits_{l=1}^M\sum\limits_{j=1}^N\sup\limits_{t\in\mathscr{R}}|a_{ijl}(t)|B_j,\sup\limits_{t\in\mathscr{R}}|J_i(t)|\Big],\nn\\
&&\frac{\epsilon}{5}\inf\limits_{t\in\mathscr{R}}c_i(t)
\Big/\sum\limits_{l=1}^M\sum\limits_{j=1}^N\sup\limits_{t\in\mathscr{R}}|a_{ijl}(t)|\sup\limits_{\zeta\in\mathscr{R}}|\dot{g}_j(\zeta)|,\,\,\,
\frac{\epsilon}{5}\inf\limits_{t\in\mathscr{R}}c_i(t)\Big/M\sum\limits_{j=1}^NB_j,\,\,\,\,
\frac{\epsilon}{5}\inf\limits_{t\in\mathscr{R}}c_i(t)\Big\}.\nn
\end{ea}
From basic properties $(\blacktriangle_2)$ and $(\blacktriangle_3)$,
we know that there exists a positive constant $l(\epsilon^*)$ such
that any interval $[\varsigma,\varsigma+l]$
$(\varsigma\in\mathscr{R})$ contains at least one common
$\epsilon^*$-almost period ${\mathcal{T}}$, namely
\begin{ea}
\left.\begin{ar}{l}
\displaystyle\Big|\phi^j\Big(t+{\mathcal{T}}-s-\kappa_{ijl}(t+{\mathcal{T}}-s)\Big)
-\phi^j\Big(t-s-\kappa_{ijl}(t-s)\Big)\Big|
_{\left.\begin{array}{l}\\\end{array}\right.}\\
\displaystyle\leq\epsilon^*<\frac{\epsilon}{5}\inf\limits_{t\in\mathscr{R}}c_i(t)
\Big/\sum\limits_{l=1}^M\sum\limits_{j=1}^N\sup\limits_{t\in\mathscr{R}}|a_{ijl}(t)|
\sup\limits_{\zeta\in\mathscr{R}}|\dot{g}_j(\zeta)|,
_{\left.\begin{array}{l}\\\end{array}\right.}\\
\displaystyle \Big|a_{ijl}(t+{\mathcal{T}}-s)-a_{ijl}(t-s)\Big|
\leq\epsilon^*<\frac{\epsilon}{5}\inf\limits_{t\in\mathscr{R}}c_i(t)\Big/\sum\limits_{j=1}^NMB_j,
_{\left.\begin{array}{l}\\\end{array}\right.}\\
\Big|J_i(t+{\mathcal{T}}-u)-J_i(t-u)\Big|
\displaystyle\leq\epsilon^*<\frac{\epsilon}{5}\inf\limits_{t\in\mathscr{R}}c_i(t),
_{\left.\begin{array}{l}\\\end{array}\right.}\\
\displaystyle \Big|c_i(t+{\mathcal{T}}-u)-c_i(t-u)\Big|
\leq\epsilon^*<\frac{\epsilon}{5}\inf\limits_{t\in\mathscr{R}}c_i(t)
\Big/\max\Big\{\sum\limits_{l=1}^M\sum\limits_{j=1}^N\sup\limits_{t\in\mathscr{R}}|a_{ijl}(t)|B_j,
\sup\limits_{t\in\mathscr{R}}|J_i(t)|\Big\}.
\end{ar}\right\}
\end{ea}\nid
For convenience, we define continuous functions
$\Xi_i(t,s)=\exp\Big(-\int\limits_0^sc_i(t-u)du\Big).$ From (6), we
get that
\begin{ea}
&&\Big|({\mathscr{F}}_i^\Sigma\phi)(t+{\mathcal{T}})-({\mathscr{F}}_i^\Sigma\phi)(t)\Big|\nn\\
&\hskip-2mm=&\hskip-2mm\displaystyle\Bigg|\int\limits_0^\infty\Bigg[\sum\limits_{l=1}^M\sum\limits_{j=1}^N
a_{ijl}(t+{\mathcal{T}}-s)g_j\Big(\sigma_j\phi^j(t+{\mathcal{T}}-s-\kappa_{ijl}(t+{\mathcal{T}}-s))\Big)
+J_i(t+{\mathcal{T}}-s)\Bigg]\nn\\
&\hskip-2mm\times&\hskip-2mm\exp\Big(-\int\limits_0^sc_i(t+{\mathcal{T}}-u)du\Big)ds-\int\limits_0^\infty\Bigg[\sum\limits_{l=1}^M\sum\limits_{j=1}^N
a_{ijl}(t-s)g_j\Big(\sigma_j\phi^j(t-s-\kappa_{ijl}(t-s))\Big)\nn\\
&\hskip-2mm+&\hskip-2mmJ_i(t-s)\Bigg]\exp\Big(-\int\limits_0^sc_i(t-u)du\Big)ds\Bigg|\nn\\
&\hskip-2mm=&\hskip-2mm\displaystyle\Bigg|\int\limits_0^\infty\sum\limits_{l=1}^M\sum\limits_{j=1}^N
\Bigg\{\Big[a_{ijl}(t+{\mathcal{T}}-s)-a_{ijl}(t-s)\Big]
g_j\Big(\sigma_j\phi^j(t+{\mathcal{T}}-s-\kappa_{ijl}(t+{\mathcal{T}}-s))\Big)\nn\\
&\hskip-2mm\times&\hskip-2mm\exp\Big(-\int\limits_0^sc_i(t+{\mathcal{T}}-u)du\Big)ds+
a_{ijl}(t-s)\Big[g_j\Big(\sigma_j\phi^j(t+{\mathcal{T}}-s-\kappa_{ijl}(t+{\mathcal{T}}-s))\Big)\nn\end{ea}
\begin{ea}
&\hskip-2mm-&\hskip-2mmg_j\Big(\sigma_j\phi^j(t-s-\kappa_{ijl}(t-s))\Big)\Big]
\exp\Big(-\int\limits_0^sc_i(t+{\mathcal{T}}-u)du\Big)ds\nn\\
&\hskip-2mm+&\hskip-2mma_{ijl}(t-s)g_j\Big(\sigma_j\phi^j(t-s-\kappa_{ijl}(t-s))\Big)\Big[
\exp\Big(-\int\limits_0^sc_i(t+{\mathcal{T}}-u)du\Big)ds\nn\\
&\hskip-2mm-&\hskip-2mm
\exp\Big(-\int\limits_0^sc_i(t-u)du\Big)ds\Big]\Bigg\}+\Bigg[J_i(t+{\mathcal{T}}-s)-J_i(t-s)\Bigg]
\exp\Big(-\int\limits_0^sc_i(t+{\mathcal{T}}-u)du\Big)ds\nn\\
&\hskip-2mm+&\hskip-2mmJ_i(t-s)\Bigg[
\exp\Big(-\int\limits_0^sc_i(t+{\mathcal{T}}-u)du\Big)ds-
\exp\Big(-\int\limits_0^sc_i(t-u)du\Big)ds\Bigg]\Bigg|\nn\\
&\hskip-2mm\leq&\hskip-2mm
\int\limits_0^\infty\sum\limits_{l=1}^M\sum\limits_{j=1}^N\Big|a_{ijl}(t+{\mathcal{T}}-s)-a_{ijl}(t-s)\Big|\Big|g_j\Big(\sigma_j\phi^j\big(t+{\mathcal{T}}-s
-\kappa_{ijl}(t+{\mathcal{T}}-s)\big)\Big)\Big|\Xi_i(t+\mathcal{T},s)ds\nn\\
&\hskip-2mm+&\hskip-2mm
\int\limits_0^{\infty}\sum\limits_{l=1}^M\sum\limits_{j=1}^N|a_{ijl}(t-s)|\Big|g_j\Big(\sigma_j\phi^j\big(t+{\mathcal{T}}-s-\kappa_{ijl}(t+{\mathcal{T}}-s)\big)\Big)
-g_j\Big(\sigma_j\phi^j\big(t-s-\kappa_{ijl}(t-s)\big)\Big)\Big|\nn\\
&\hskip-2mm\times&\hskip-2mm \Xi_i(t+\mathcal{T},s)ds+
\int\limits_0^{\infty}\sum\limits_{l=1}^M\sum\limits_{j=1}^N|a_{ijl}(t-s)|\Big|g_j\Big(\sigma_j\phi^j\big(t-s-\kappa_{ijl}(t-s)\big)\Big)\Big|
\Big|\Xi_i(t+\mathcal{T},s) -\Xi_i(t,s)\Big|ds\nn\\
&\hskip-2mm+&\hskip-2mm
\int\limits_0^{\infty}\Big|J_i(t-s)\Big|\Big|\Xi_i(t+\mathcal{T},s)
-\Xi_i(t,s)\Big|ds+
\int\limits_0^{\infty}\Big|J_i(t+{\mathcal{T}}-s)-J_i(t-s)\Big|\Xi_i(t+\mathcal{T},s)ds.
\end{ea}\nid
By using $(H_1)$, (9)-(10) and mean value theorem of differential
calculus, we obtain that
\begin{ea}
&&\Big|({\mathscr{F}}_i^\Sigma\phi)(t+{\mathcal{T}})-({\mathscr{F}}_i^\Sigma\phi)(t)\Big|
\leq\epsilon^*\sum\limits_{l=1}^M\sum\limits_{j=1}^NB_j\int\limits_0^\infty
\exp\Big(-\int\limits_0^s\inf\limits_{t\in\mathscr{R}}c_i(t)du\Big)ds+\sum\limits_{l=1}^M\sum\limits_{j=1}^N\sup\limits_{t\in\mathscr{R}}|a_{ijl}(t)|\nn\\
&\hskip-2mm\times&\hskip-2mm
\int\limits_0^{\infty}|\dot{g}_j(\theta)|
\sigma_j\Big|\phi^j\big(t+{\mathcal{T}}-s-\kappa_{ijl}(t+{\mathcal{T}}-s)\big)
-\phi^j\big(t-s-\kappa_{ijl}(t-s)\big)\Big|\exp\Big(-\int\limits_0^s
\inf\limits_{t\in\mathscr{R}}c_i(t)du\Big)ds\nn\\
&\hskip-2mm+&\hskip-2mm
\sum\limits_{l=1}^M\sum\limits_{j=1}^N\Big(\sup\limits_{t\in\mathscr{R}}|a_{ijl}(t)|B_j+\sup\limits_{t\in\mathscr{R}}|J_i(t)|\Big)\int\limits_0^{\infty}
\int\limits_0^s\exp(\widetilde{\theta})\Big|c_i(t+{\mathcal{T}}-u)-c_i(t-u)\Big|du
ds\nn\\&\hskip-2mm+&\hskip-2mm
\epsilon^*\int\limits_0^{\infty}\exp\Big(-\int\limits_0^s\inf\limits_{t\in\mathscr{R}}c_i(t)du\Big)ds,
\end{ea}\nid
where $\theta$ lies between
$\sigma_j\phi^j\big(t+{\mathcal{T}}-s-\kappa_{ijl}(t+{\mathcal{T}}-s)\big)$
and $\sigma_j\phi^j\big(t-s-\kappa_{ijl}(t-s)\big)$,
$\widetilde{\theta}$ lies between
$-\int\limits_0^sc_i(t+{\mathcal{T}}-u)du$ and
$-\int\limits_0^sc_i(t-u)du$. Noting that
$-\int\limits_0^s\sup\limits_{t\in\mathscr{R}}c_i(t)du\leq\widetilde{\theta}\leq
-\int\limits_0^s\inf\limits_{t\in\mathscr{R}}c_i(t)du,$ it follows
that
$\exp(\widetilde{\theta})\leq\exp(-s\inf\limits_{t\in\mathscr{R}}c_i(t)).$
Therefore, by using (9) and (11), we have
\begin{ea}
\Big|({\mathscr{F}}_i^\Sigma\phi)(t+{\mathcal{T}})-({\mathscr{F}}_i^\Sigma\phi)(t)\Big|
&\hskip-2mm\leq&\hskip-2mm\epsilon^*\sum\limits_{j=1}^NMB_j\Big/\inf\limits_{t\in\mathscr{R}}c_i(t)+
\epsilon^*\int\limits_0^{\infty}\exp\Big(-\int\limits_0^s\inf\limits_{t\in\mathscr{R}}c_i(t)du\Big)ds
\nn\\
&\hskip-2mm+&\hskip-2mm
\epsilon^*\sum\limits_{l=1}^M\sum\limits_{j=1}^N(\sup\limits_{t\in\mathscr{R}}|a_{ijl}(t)|B_j+\sup\limits_{t\in\mathscr{R}}|J_i(t)|)\int\limits_0^{\infty}
\exp(-s\inf\limits_{t\in\mathscr{R}}c_i(t))sds\nn\\
&\hskip-2mm+&\hskip-2mm
\epsilon^*\sum\limits_{l=1}^M\sum\limits_{j=1}^N\sigma_j\sup\limits_{t\in\mathscr{R}}|a_{ijl}(t)|\sup\limits_{\zeta\in\mathscr{R}}
|\dot{g}_j(\zeta)|\Big/\inf\limits_{t\in\mathscr{R}}c_i(t)\leq\epsilon,
\end{ea}\nid
which leads to the almost periodicity of $\mathscr{F}^\Sigma(\phi)$.
From Step 1, it follows that
${\mathscr{F}^\Sigma}:{\mathscr{B}}^\Sigma\rightarrow{\mathscr{B}}^\Sigma$.
That is, each ${\mathscr{B}}^\Sigma$ is invariant region of
$\mathscr{F}^\Sigma$.

{\bf Step 3:}\,\,We should prove that
${\mathscr{F}^\Sigma}:{\mathscr{B}}^\Sigma\rightarrow{\mathscr{B}}^\Sigma$
is continuous. Take any two $\phi_1,\phi_2\in{\mathscr{B}}^\Sigma$.
From (6) and Lagrange's mean value theorem, we have
\begin{ea}
&&\hskip-2mm\Big|({\mathscr{F}}_i^\Sigma\phi_1)(t)-({\mathscr{F}}_i^\Sigma\phi_2)(t)\Big|
\leq\int\limits_0^{\infty}\sum\limits_{l=1}^M\sum\limits_{j=1}^N|a_{ijl}(t-s)|\nn\\
&\hskip-2mm\times&\hskip-2mm\Big|g_j(\sigma_j\phi_1^j(t-s-\kappa_{ijl}(t-s)))-g_j(\sigma_j\phi_2^j(t-s-\kappa_{ijl}(t-s)))\Big|
\exp\Big(-\int\limits_0^sc_i(t-u)du\Big)ds\nn\\
&\hskip-2mm\leq&\hskip-2mm\sum\limits_{l=1}^M\sum\limits_{j=1}^N\sup\limits_{t\in\mathscr{R}}|a_{ijl}(t)|
\sigma_j\sup\limits_{\zeta\in\mathscr{R}}\dot{g}_j(\zeta)\Big/
\inf\limits_{t\in\mathscr{R}}c_i(t)\|\phi_1-\phi_2\|,\nn
\end{ea}
which leads to
$$
\|{\mathscr{F}^\Sigma}\phi_1-{\mathscr{F}^\Sigma}\phi_2\|
\leq\max\limits_{i\in\mathscr{I}}\Big\{\sum\limits_{l=1}^M\sum\limits_{j=1}^N\sup\limits_{t\in\mathscr{R}}|a_{ijl}(t)|
\sigma_j\sup\limits_{\zeta\in\mathscr{R}}\dot{g}_j(\zeta)\Big/
\inf\limits_{t\in\mathscr{R}}c_i(t)\Big\}\|\phi_1-\phi_2\|.
$$
This implies that ${\mathscr{F}^\Sigma}$ is continuous with respect
to $\phi\in{\mathscr{B}}^\Sigma$.

From Lemma 3, each ${\mathscr{B}}^\Sigma$ is compact convex subset.
Since
${\mathscr{F}^\Sigma}:{\mathscr{B}}^\Sigma\rightarrow{\mathscr{B}}^\Sigma$
is continuous, by Schauder's fixed point theorem, there exists at
least one $\widehat{u}_\Sigma\in{\mathscr{B}}^\Sigma$ such that
${\mathscr{F}^\Sigma}\widehat{u}_\Sigma=\widehat{u}_\Sigma$. It is
easy for us to check that
\begin{ea}
\frac{d\widehat{u}^i_\Sigma(t)}{dt}
&\hskip-2mm=&\hskip-2mm\frac{d}{dt}\displaystyle\int\limits_{-\infty}^t\Big[
\sum\limits_{l=1}^M\sum\limits_{j=1}^Na_{ijl}(s)g_j\Big(\sigma_j\widehat{u}^j_\Sigma(s-\kappa_{ijl}(s))\Big)+J_i(s)\Big]
\exp\Big(-\int\limits_s^tc_i(u)du\Big)ds\nn\\
&=&\hskip-2mm-c_i(t)\widehat{u}^i_\Sigma(t)+\sum\limits_{l=1}^M\sum\limits_{j=1}^Na_{ijl}(t)
g_j\Big(\sigma_j\widehat{u}^j_\Sigma(t-\kappa_{ijl}(t))\Big)+J_i(t).
\end{ea}
Hence $\widehat{u}_\Sigma(t)$ is an almost periodic solution of GNNs
(1) in ${\mathscr{B}}^{\Sigma}$. The proof is complete. $\square$

\vskip 0.1in From above theorem, we know that each
${\mathscr{B}}^{\Sigma}$ is an invariant region of GNNs (1) and
there exist at least $2^N$ almost periodic encoded patterns in these
${\mathscr{B}}^{\Sigma}$. In what follows, we should prove that each
$\mathscr{K}^\Sigma$ is an attracting basin for almost periodic
encoded pattern in ${\mathscr{B}}^{\Sigma}$ under our additional
assumptions. For convenience, we take the following notations:
$$\mathcal{C}=diag\Big(\inf\limits_{t\in\mathscr{R}}c_1(t),
\inf\limits_{t\in\mathscr{R}}c_2(t),\cdots,\inf\limits_{t\in\mathscr{R}}c_N(t)\Big),\,\,\,\,\,
\mathcal{H}=\Big(\hbar_{ij}\Big)_{N\times
N}\,\,\,\,\mbox{with}\,\,\,\hbar_{ij}:=\sum\limits_{l=1}^M\sup\limits_{t\in\mathscr{R}}|a_{ijl}(t)|,$$
$$\mathcal{G}=diag\Big(\sigma_1\dot{g}_1(\zeta),
\sigma_2\dot{g}_2(\zeta),\cdots,\sigma_N\dot{g}_N(\zeta)\Big)\,\,\,
\mbox{with}\,\,\,
\dot{g}_i(\zeta)=\max\Big\{\dot{g}_i(z)\Big|z=\sigma_i\beta_{i1},\sigma_i\alpha_{i2}\Big\}.$$
Next, we should introduce two additional assumptions: \vskip 0.1in
$\bullet$$(H_3^\mathcal{A})$\ \
$\mathcal{C}-\mathcal{H}\mathcal{G}\in\mathscr{M}$. \vskip 0.1in
$\bullet$$(H_4^\mathcal{A})$\ \ There exist positive constants
$p>1$, $d_i$, $o_k$, $q_{jlk}$ and $p_{jlk}$ such that
\begin{ea*}
pd_i\inf\limits_{t\in\mathscr{R}}c_i(t)
\hskip-2mm&>&\hskip-2mm\ss\sum\limits_{l=1}^M\sum\limits_{j=1}^N\Big[
d_j\big(\sup\limits_{t\in\mathscr{R}}|a_{jil}(t)|\big)^{pp_{il,m+1}}\big(\sigma_i\dot{g}_i(\zeta)\big)^{pq_{il,m+1}}\nn\\
\hskip-2mm&+&\hskip-2mm\ss\sum\limits_{k=1}^md_io_k\big(\sup\limits_{t\in\mathscr{R}}|a_{ijl}(t)|\big)^{pp_{jlk}/o_k}
\big(\sigma_j\dot{g}_j(\zeta)\big)^{pq_{jlk}/o_k}\Big],\nn
\end{ea*}
where
$\dot{g}_i(\zeta)=\max\Big\{\dot{g}_i(z)\Big|z=\sigma_i\beta_{i1},\sigma_i\alpha_{i2}\Big\}$,
$\sum\limits_{k=1}^mo_k=p-1$,
$\sum\limits_{k=1}^{m+1}q_{jlk}=\sum\limits_{k=1}^{m+1}p_{jlk}=1$
for each $j\in\mathscr{I}$ and $l\in\mathscr{L}$, $m$ is a positive
integer.

\vskip 0.1in\noindent\textsl{Theorem 3}\\Under the basic assumptions
$(H_1)$ and $(H_1^\mathcal{A})$-$(H_3^\mathcal{A})$, for each
$\Sigma$, there exists a unique almost periodic encoded pattern of
GNNs (1) which is exponentially stable in $\mathscr{K}^\Sigma$.

\vskip 0.1in\noindent \textsl{Proof}\\From Theorem 2, there exists
 an almost periodic solution $\widehat{u}$ of GNNs (1) in each
${\mathscr{B}}^\Sigma$. For any initial condition
$\phi\in\mathscr{K}^\Sigma$, by Theorem 1, we know that
$x_t(\cdot;\phi)\in\mathscr{K}^\Sigma$ for all $t\geq0$. Under
translation $y(t)=\widehat{u}(t)-x(t;\phi)$, we get that
\begin{ea}
\frac{dy^i(t)}{dt}=-c_i(t)y^i(t)+\sum\limits_{l=1}^M\sum\limits_{j=1}^Na_{ijl}(t)
\Big[g_j\big(\sigma_j\widehat{u}^j(t-\kappa_{ijl}(t))\big)-g_j\big(\sigma_jx^j(t-\kappa_{ijl}(t))\big)\Big].
\end{ea}
From $(H_3^\mathcal{A})$ and Lemma 5, there exists a
$K=(K_1,K_2,\cdots,K_N)^T\in\Omega_\mathscr{M}(\mathcal{C}-\mathcal{H}\mathcal{G})$
such that
$\sup\limits_{\theta\in[-\kappa,0]}|\widehat{u}^i(\theta)-\phi^i(\theta)|\leq
K_i$. Then we get that
\begin{ea*}\inf\limits_{t\in\mathscr{R}}c_i(t)K_i-\sum\limits_{l=1}^M\sum\limits_{j=1}^N
\sup\limits_{t\in\mathscr{R}}|a_{ijl}(t)|\sigma_j\dot{g}_j(\zeta)K_j>0,
\end{ea*}where
$\dot{g}_j(\zeta)=\max\Big\{\dot{g}_j(z)\Big|z=\sigma_j\beta_{j1},\sigma_j\alpha_{j2}\Big\}$.
We consider the single-variable functions $W_i(\cdot)$ defined by
$$W_i(\theta)=(\inf\limits_{t\in\mathscr{R}}c_i(t)-\theta)K_i-\sum\limits_{l=1}^M\sum\limits_{j=1}^N\sup\limits_{t\in\mathscr{R}}|a_{ijl}(t)|
\sigma_j\dot{g}_j(\zeta)K_je^{\theta\kappa_{ijl}}.$$ Noting that
$W_i(0)>0$ and $W_i(\theta)\rightarrow-\infty$ as
$\theta\rightarrow+\infty$, there exists a suitable $\mu$ such that
for all $i\in\mathscr{I}$,
\begin{ea}(\inf\limits_{t\in\mathscr{R}}c_i(t)-\mu)K_i-\sum\limits_{l=1}^M\sum\limits_{j=1}^N\sup\limits_{t\in\mathscr{R}}|a_{ijl}(t)|
\sigma_j\dot{g}_j(\zeta)K_je^{\mu\kappa_{ijl}}>0.
\end{ea}
\indent Consider function $Z_i(t)=e^{{\mu}t}|y^i(t)|$ where
$t\in[-\kappa,\infty)$. Let $\varrho>1$. It is obvious that
$Z_i(t)<\varrho K_i$ for all $t\in[-\kappa,0]$. Now we claim that
$Z_i(t)<\varrho K_i$ for all $t>0$ and $i\in\mathscr{I}$. Otherwise
there is a first time $t_0>0$ and some $i^*\in\mathscr{I}$ such that
$Z_{i^*}(t_0)=\varrho K_{i^*}$, $\frac{d^+|Z_{i^*}(t_0)|}{dt}\geq0$
and $Z_j(t)<\varrho K_j$ ($j\neq {i^*}$) for all
$t\in[-\kappa,t_0]$. From (14), we derive that
\begin{ea}
\frac{d^+|y^{i^*}(t_0)|}{dt}\leq-\inf\limits_{t\in\mathscr{R}}c_{i^*}(t)|y^{i^*}(t_0)|
+\sum\limits_{l=1}^M\sum\limits_{j=1}^N\sup\limits_{t\in\mathscr{R}}|a_{{i^*}jl}(t)|\sigma_j\dot{g}_j(\xi)\big|y^j(t_0-\kappa_{{i^*}jl}(t_0))\big|,\nn
\end{ea}
where $\xi$ lies between
$\sigma_j\widehat{u}^j(t_0-\kappa_{{i^*}jl}(t_0))$ and
$\sigma_jx^j(t_0-\kappa_{{i^*}jl}(t_0))$. From above inequality and
(15), we get that
\begin{ea*}
\frac{d^+|Z_{i^*}(t_0)|}{dt}
\hskip-2mm&\leq&\hskip-2mm-\big(\inf\limits_{t\in\mathscr{R}}c_{i^*}(t)-\mu\big)Z_{i^*}(t_0)
+\sum\limits_{l=1}^M\sum\limits_{j=1}^N\sup\limits_{t\in\mathscr{R}}|a_{{i^*}jl}(t)|\sigma_j\dot{g}_j(\zeta)
e^{\mu\kappa_{{i^*}jl}(t_0)}Z_j(t_0-\kappa_{{i^*}jl}(t_0))\nn\\
\hskip-2mm&\leq&\hskip-2mm-\big(\inf\limits_{t\in\mathscr{R}}c_{i^*}(t)-\mu\big)Z_{i^*}(t_0)
+\sum\limits_{l=1}^M\sum\limits_{j=1}^N\sup\limits_{t\in\mathscr{R}}|a_{{i^*}jl}(t)|\sigma_j\dot{g}_j(\zeta)
e^{\mu\kappa_{{i^*}jl}}\sup\limits_{\theta\in[t_0-\kappa,t_0]}Z_j(\theta)\nn\\
\hskip-2mm&\leq&\hskip-2mm
-\big(\inf\limits_{t\in\mathscr{R}}c_{i^*}(t)-\mu\big)\varrho
K_{i^*}+\sum\limits_{l=1}^M\sum\limits_{j=1}^N\sup\limits_{t\in\mathscr{R}}|a_{{i^*}jl}(t)|
\sigma_j\dot{g}_j(\zeta)e^{\mu\kappa_{{i^*}jl}}\varrho K_j<0,
\end{ea*}
which leads to a contradiction. Hence $Z_i(t)<\varrho K_i$ for all
$t>0$ and $i\in\mathscr{I}$. That is, there exists a positive
constant $\widetilde{\varrho}$ such that
$$|\widehat{u}^i(t)-x^i(t;\phi)|\leq e^{-\mu t}\widetilde{\varrho}\sup\limits_{\theta\in[-\kappa,0]}|\widehat{u}^i(\theta)-\phi^i(\theta)|,$$
for all $t\geq0$ and $i\in\mathscr{I}$. Therefore, for each
$\Sigma$, there exists a unique almost periodic encoded pattern
$\widehat{u}(t)$ which is exponentially stable  in
${\mathscr{K}}^\Sigma$. The proof is complete.
\indent\indent\indent\indent$\square$
%\vskip 0.1in\noindent{\bf Remark 3.1.} For activation functions of
%class $\mathcal{A}$, we only assume that $|g_i(x)-g_i(y)|\leq
%\widetilde{L}_i|x-y|$ and $|g_i(x)|\leq \widetilde{B}_i$ for all
%$x,y\in\mathscr{R}$, where $\widetilde{L}_i$ and $\widetilde{B}_i$
%are positive constants. Under the assumption $(H_3^\mathcal{A})$
%with $\mathcal{G}=diag\Big(\sigma_1\widetilde{L}_1,
%\sigma_2\widetilde{L}_2,\cdots,\sigma_N\widetilde{L}_N\Big)$,
%similarly as Lemma 2.1, Theorem 3.1 and Theorem 3.2, there exists a
%unique stable almost periodic encoded pattern of GNNs (1.1) which
%includes the related results in [10].

\vskip 0.1in\noindent\textsl{Theorem 4}\\Under the basic assumptions
$(H_1)$, $(H_1^\mathcal{A})$-$(H_2^\mathcal{A})$ and
$(H_4^\mathcal{A})$, for each $\Sigma$, there exists a unique almost
periodic encoded pattern of GNNs (1) which is exponentially stable
in $\mathscr{K}^\Sigma$.

\vskip 0.1in\noindent \textsl{Proof}\\From Theorem 2, there exists
 an almost periodic solution $\widehat{u}$ of GNNs (1) in each
${\mathscr{B}}^\Sigma$. For any $\phi\in\mathscr{K}^\Sigma$, by
Theorem 1, we know that $x_t(\cdot;\phi)\in\mathscr{K}^\Sigma$ for
all $t\geq0$. Let $y(t)=\widehat{u}(t)-x(t;\phi)$. We consider the
Lyapunov functional $V(y)(t)=\sum\limits_{i=1}^Nd_i|y^i(t)|^p$,
where $d_i>0$. From (14), we can derive that
\begin{ea}
&&\hskip-6mm\frac{d^+V(y)(t)}{dt}\leq\sum\limits_{i=1}^Npd_i|y^i(t)|^{p-1}\big\{
-\inf\limits_{t\in\mathscr{R}}c_i(t)|y^i(t)|
+\sum\limits_{l=1}^M\sum\limits_{j=1}^N\sup\limits_{t\in\mathscr{R}}|a_{ijl}(t)|
\sigma_j\dot{g}_j(\zeta)\big|y^j(t-\kappa_{ijl}(t))\big|\big\}\nn\\
&\hskip-2mm\leq&\hskip-2mm-p\sum\limits_{i=1}^N\Big\{\inf\limits_{t\in\mathscr{R}}c_i(t)d_i|y^i(t)|^{p}
-\sum\limits_{l=1}^M\sum\limits_{j=1}^Nd_i|y^i(t)|^{p-1}\sup\limits_{t\in\mathscr{R}}|a_{ijl}(t)|
\sigma_j\dot{g}_j(\zeta)\big|y^j(t-\kappa_{ijl}(t))\big|\Big\}\nn\\
&\hskip-2mm\leq&\hskip-2mm-p\sum\limits_{i=1}^N\Big\{\inf\limits_{t\in\mathscr{R}}c_i(t)d_i|y^i(t)|^{p}
-\sum\limits_{l=1}^M\sum\limits_{j=1}^Nd_i\prod\limits_{k=1}^m
\Big[\big(\sup\limits_{t\in\mathscr{R}}|a_{ijl}(t)|\big)^{p_{jlk}/o_k}
\big(\sigma_j\dot{g}_j(\zeta)\big)^{q_{jlk}/o_k}|y^i(t)|\Big]^{o_k}\nn\\
&\hskip-2mm\times&\hskip-2mm
\Big[\big(\sup\limits_{t\in\mathscr{R}}|a_{ijl}(t)|\big)^{p_{jl,m+1}}
\big(\sigma_j\dot{g}_j(\zeta)\big)^{q_{jl,m+1}}\big|y^j(t-\kappa_{ijl}(t))\big|\Big]
\Big\}\nn\\
&\hskip-2mm\leq&\hskip-2mm-p\sum\limits_{i=1}^N\Big\{\inf\limits_{t\in\mathscr{R}}c_i(t)d_i|y^i(t)|^{p}
-\sum\limits_{l=1}^M\sum\limits_{j=1}^N
\Big[\frac{d_i}{p}\sum\limits_{k=1}^mo_k\big(\sup\limits_{t\in\mathscr{R}}|a_{ijl}(t)|\big)^{pp_{jlk}/o_k}
\big(\sigma_j\dot{g}_j(\zeta)\big)^{pq_{jlk}/o_k}|y^i(t)|^p\nn\\
&\hskip-2mm+&\hskip-2mm
\frac{d_i}{p}\big(\sup\limits_{t\in\mathscr{R}}|a_{ijl}(t)|\big)^{pp_{jl,m+1}}
\big(\sigma_j\dot{g}_j(\zeta)\big)^{pq_{jl,m+1}}\big|y^j(t-\kappa_{ijl}(t))\big|^p\Big]
\Big\}\nn\\
&\hskip-2mm\leq&\hskip-2mm-\sum\limits_{i=1}^N\Big[pd_i\inf\limits_{t\in\mathscr{R}}c_i(t)
-\sum\limits_{l=1}^M\sum\limits_{j=1}^N
\sum\limits_{k=1}^md_io_k\big(\sup\limits_{t\in\mathscr{R}}|a_{ijl}(t)|\big)^{pp_{jlk}/o_k}
\big(\sigma_j\dot{g}_j(\zeta)\big)^{pq_{jlk}/o_k}\Big]
|y^i(t)|^{p}\nn\\
&\hskip-2mm+&\hskip-2mm
\sum\limits_{i=1}^N\Big[\sum\limits_{l=1}^M\sum\limits_{j=1}^Nd_j
\big(\sup\limits_{t\in\mathscr{R}}|a_{jil}(t)|\big)^{pp_{il,m+1}}
\big(\sigma_i\dot{g}_i(\zeta)\big)^{pq_{il,m+1}}\Big]
\big|y^i(t-\kappa_{jil}(t))\big|^p\nn\\
&\hskip-2mm\leq&\hskip-2mm-{\alpha}V(y)(t)+{\beta}\sup\limits_{t-\kappa{\leq}s{\leq}t}V(y)(s),
\end{ea}
where
$$\alpha=\min\Big\{pd_i\inf\limits_{t\in\mathscr{R}}c_i(t)
-\sum\limits_{l=1}^M\sum\limits_{j=1}^N
\sum\limits_{k=1}^md_io_k\big(\sup\limits_{t\in\mathscr{R}}|a_{ijl}(t)|\big)^{pp_{jlk}/o_k}
\big(\sigma_j\dot{g}_j(\xi)\big)^{pq_{jlk}/o_k}\Big\},$$
$$\beta=\max\Big\{\sum\limits_{l=1}^M\sum\limits_{j=1}^Nd_j
\big(\sup\limits_{t\in\mathscr{R}}|a_{jil}(t)|\big)^{pp_{il,m+1}}
\big(\sigma_i\dot{g}_i(\xi)\big)^{pq_{il,m+1}}\Big\}.$$ From
$(H_4^\mathcal{A})$, we have $\alpha>\beta>0$. By using Halanay
inequality, we get for all $t\in\mathscr{R}$,
\begin{ea}
V(y)(t)\hskip-2mm&\leq&\hskip-2mm\Big(\sup\limits_{-\kappa{\leq}s{\leq}0}V(y)(s)\Big)\exp(-{\gamma}t),
\end{ea}
where $\gamma=\alpha-\beta e^{\gamma\kappa}$. It follows that
$$\sum\limits_{i=1}^N\Big|\widehat{u}^i(t)-x^i(t;\phi)\Big|^p\leq e^{-\gamma t}\frac{\max\{d_i\}}{\min\{d_i\}}
\sum\limits_{i=1}^N\sup\limits_{\theta\in[-\kappa,0]}\Big|\widehat{u}^i(\theta)-\phi^i(\theta)\Big|^p.$$
Hence $\widehat{u}(t)$ is exponentially stable. The proof is
complete.\indent\indent\indent\indent\indent\indent\indent\indent\indent$\square$

\vskip 0.1in\noindent\textsl{Remark 2}\\When $c_i$, $a_{ijl}$,
$\kappa_{ijl}$, $J_i:\mathscr{R}\rightarrow\mathscr{R}$ are
$\omega$-periodic functions with $\omega>0$, we also obtain the
existence and exponential stability of $2^N$ periodic solutions of
GNNs (1). Our results in Theorem 3 and Theorem 4 are distinguished
from the existing results on the following points: (i)Most of the
previous results of neural networks only focus on the existence and
stability of unique almost periodic (periodic) solution.  Hence, we
extend the related results [20-31,37,39,41] to the convergence
analysis of multiple almost periodic (periodic) solutions. (ii)We
not only establish existing regions for almost periodic solutions,
but also estimate attracting basins of these almost periodic
solutions. (iii)Our sufficient conditions
$(H_3^\mathcal{A})$-$(H_4^\mathcal{A})$ are dependent of system
parameters and derivative of activation functions on boundary points
which make our results new in the literature.

\vskip 0.1in\noindent\textsl{Remark 3}\\From $(H_2^\mathcal{A})$, we
know that $
F_i(z_{i2})-\sum\limits_{l=1}^M\sum\limits_{j\neq{i}}\sup\limits_{t\in\mathscr{R}}|a_{ijl}(t)|B_j
+\inf\limits_{t\in\mathscr{R}}J_i(t)>0. $ It is obvious that
$F_i(z_{i2})+\sum\limits_{l=1}^M\sum\limits_{j\neq{i}}\sup\limits_{t\in\mathscr{R}}|a_{ijl}(t)|B_j
+\sup\limits_{t\in\mathscr{R}}J_i(t)>0$. Together with (2), there
exists a $\widetilde{z}_{i1}$ with $\widetilde{z}_{i1}>z_{i1}$ such
that
$F_i(\widetilde{z}_{i1})+\sum\limits_{l=1}^M\sum\limits_{j\neq{i}}\sup\limits_{t\in\mathscr{R}}|a_{ijl}(t)|B_j
+\sup\limits_{t\in\mathscr{R}}J_i(t)=0$. Similarly, there exists a
$\widehat{z}_{i2}$ with $\widehat{z}_{i2}<z_{i2}$ such that
$F_i(\widehat{z}_{i2})-\sum\limits_{l=1}^M\sum\limits_{j\neq{i}}\sup\limits_{t\in\mathscr{R}}|a_{ijl}(t)|B_j
+\inf\limits_{t\in\mathscr{R}}J_i(t)=0$. For each
$\Sigma=(\varsigma_1,\varsigma_2,\cdots,\varsigma_N)$, we denote
$$\Xi^\Sigma:=\Big\{z=(z_1,z_2,\cdots,z_N)^T\Big|
z_i\in[\widehat{z}_{i\varsigma_i},\widetilde{z}_{i\varsigma_i}],\,\,i=1,2,\cdots,N
\,\,\, \mbox{and}\,\,\,
\mathcal{C}-\mathcal{H}\mathcal{G}\in\mathscr{M}\Big\},$$ where
$\mathcal{G}=diag\Big(\sigma_1\dot{g}_1(\sigma_1z_1),
\sigma_2\dot{g}_2(\sigma_2z_2),\cdots,\sigma_N\dot{g}_N(\sigma_Nz_N)\Big)$.
Under the basic assumptions of Theorem 3, we have
$\Xi^\Sigma\not=\emptyset$. Assume that there exists a
$\overline{z}\in\Xi^\Sigma$ such that (i)$\varsigma_i=1$,
$\overline{z}_{i}\geq\sup\limits_{z\in\Xi^\Sigma}\Big\{z_i\Big\}$
($\overline{z}_{i}\geq \widehat{z}_{i1}$); (ii)$\varsigma_i=2$,
$\overline{z}_{i}\leq\inf\limits_{z\in\Xi^\Sigma}\Big\{z_i\Big\}$
($\overline{z}_{i}\leq\widetilde{z}_{i2}$). For each $\Sigma$, we
let (i)$\varsigma_i=1$, $\beta_{i1}:=\overline{z}_{i}$;
(ii)$\varsigma_i=2$, $\alpha_{i2}:=\overline{z}_{i}$. From Theorem 1
and Theorem 3, we know that each $
\mathscr{K}^\Sigma:=\underbrace{\mathscr{K}_{1\varsigma_1}\times\mathscr{K}_{2\varsigma_2}
\cdots\times\mathscr{K}_{N\varsigma_N}}_N$ is a {\bf larger
attracting basin}.

\vskip 0.3in \normalsize \baselineskip16pt
\centerline{4.~~\textrm{SOME GENERALIZATIONS AND IMPROVEMENTS}}
\vskip 0.1in\noindent In this section, we should make some
generalizations  by considering the second class of saturated
activation functions and make some improvements by applying our
results to some special cases. The second class of activation
functions we considered in this paper satisfies
\begin{ea*}
\mbox{Class}\,\mathcal{B}\,:\,\,\, g_j\in\mathcal{C},\,\,\,
g_j(x)=\left\{\begin{aligned}
u_{j1}& \,\,\,\,\textrm{if}\,\,\,\, -\infty<x<\ell_{j1}, \\
\widetilde{g}_j(x)& \,\,\,\,\mbox{if}\,\,\,\,  \ell_{j1}\leq{x}\leq\ell_{j2},\\
u_{j2}& \,\,\,\,\mbox{if}\,\,\,\,  \ell_{j2}<x<+\infty,
\end{aligned} \right.
\end{ea*}
where $\widetilde{g}_j\in\mathcal{C}^1$ is an increasing function
with $\widetilde{g}_j(0)=0$, $\ell_{j1}<0<\ell_{j2}$ and
$-\infty<u_{j1}<0<u_{j2}<+\infty$. Similarly as Lemma 1, it is easy
for us to have the following lemma.

\vskip 0.1in\noindent\textsl{Lemma 6}\\Assume that $(H_1)$ holds.
For any
$\phi=(\phi^1,\phi^2,\cdots,\phi^N)^T\in{\mathcal{C}}([-\kappa,0],\mathscr{R}^N)$,
$$\|\phi^i\|_\kappa\leq\Big(\sum\limits_{l=1}^M\sum\limits_{j=1}^N\sup\limits_{t\in\mathscr{R}}|a_{ijl}(t)|
\cdot\max\limits_{k=1,2}\Big\{|u_{jk}|\Big\}
+\sup\limits_{t\in\mathscr{R}}|J_i(t)|\Big)
\Big/\inf\limits_{t\in\mathscr{R}}c_i(t)$$ implies that
$$\|u_t^i(\cdot;\phi)\|_\kappa\leq\Big(\sum\limits_{l=1}^M\sum\limits_{j=1}^N\sup\limits_{t\in\mathscr{R}}|a_{ijl}(t)|
\cdot\max\limits_{k=1,2}\Big\{|u_{jk}|\Big\}
+\sup\limits_{t\in\mathscr{R}}|J_i(t)|\Big)
\Big/\inf\limits_{t\in\mathscr{R}}c_i(t)$$ for all $t\geq 0$, where
$u(t;\phi)$ is the solution of GNNs (1) with $u_0(s)=\phi(s)$ for
$s\in[-\kappa,0]$.

\vskip 0.1in\noindent\textsl{Remark 4}\\Since each activation
function of Class $\mathcal{B}$ is bounded with
$|g_j(x)|\leq\max\limits_{k=1,2}\Big\{|u_{jk}|\Big\}$,
$j\in\mathscr{I}$. As similar proof of Lemma 1, it is easy for us to
get the result of Lemma 6.

\vskip 0.1in\indent For activation functions of class $\mathcal{B}$,
we consider the following two parameter assumptions which are used
to establish the existence and stability of $2^N$ almost periodic
encoded patterns of GNNs (1): \vskip 0.1in
$\bullet(H_1^\mathcal{B}):$$\,\,\displaystyle
\sup\limits_{t\in\mathscr{R}}c_i(t)<\inf\limits_{t\in\mathscr{R}}
\sum\limits_{l=1}^Ma_{iil}(t)\sigma_i\dot{\widetilde{g}}_i(\zeta),\,\,\,\zeta\in[\ell_{i1},\ell_{i2}].$

\vskip 0.1in
$\bullet(H_2^\mathcal{B}):$$\ss(-1)^k\cdot\Big\{-\sup\limits_{t\in\mathscr{R}}c_i(t)\frac{\ell_{ik}}{\sigma_i}+\inf\limits_{t\in\mathscr{R}}\sum\limits_{l=1}^Ma_{iil}(t)u_{ik}
+J_i(t)\Big\}$\\
\indent\indent\indent$\ss>\sum\limits_{l=1}^M\sum\limits_{j\neq{i}}\sup\limits_{t\in\mathscr{R}}|a_{ijl}(t)|\cdot\max\limits_{k=1,2}\Big\{|u_{jk}|\Big\}$\\
for all $t\in\mathscr{R}$, where $i\in\mathscr{I}$ and $k=1,2$.

\vskip 0.1in \noindent Take $k=1$ in $(H_2^\mathcal{B})$, it is easy
for us to derive that
\begin{ea}
F_i(\frac{\ell_{i1}}{\sigma_i})+\sum\limits_{l=1}^M\sum\limits_{j\neq{i}}\sup\limits_{t\in\mathscr{R}}|a_{ijl}(t)|
\cdot\max\limits_{k=1,2}\Big\{|u_{jk}|\Big\}+\sup\limits_{t\in\mathscr{R}}J_i(t)<0.\nn
\end{ea}
Noting that $F_i(z)\rightarrow+\infty$ as $z\rightarrow-\infty$, we
know that there exists a $\widehat{\ell}_{i1}$ with
$\ss\widehat{\ell}_{i1}<\frac{\ell_{i1}}{\sigma_i}<0$ such that
$$F_i(\widehat{\ell}_{i1})+\sum\limits_{l=1}^M\sum\limits_{j\neq{i}}\sup\limits_{t\in\mathscr{R}}|a_{ijl}(t)|
\cdot\max\limits_{k=1,2}\Big\{|u_{jk}|\Big\}+\sup\limits_{t\in\mathscr{R}}J_i(t)=0.$$
Take $k=2$ in $(H_2^\mathcal{B})$, by the similar argument, there
exists a $\widetilde{\ell}_{i2}$ with
$\ss0<\frac{\ell_{i2}}{\sigma_i}<\widetilde{\ell}_{i2}$ such that
$$F_i(\widetilde{\ell}_{i2})-\sum\limits_{l=1}^M\sum\limits_{j\neq{i}}\sup\limits_{t\in\mathscr{R}}|a_{ijl}(t)|
\cdot\max\limits_{k=1,2}\Big\{|u_{jk}|\Big\}+\inf\limits_{t\in\mathscr{R}}J_i(t)=0.$$

\indent For convenience, as Section 3, we  also take the following
denotations:
\begin{eq}
\left\{\begin{ar}{l}\ss\alpha_{i1}=-\Big(\sum\limits_{l=1}^M\sum\limits_{j=1}^N\sup\limits_{t\in\mathscr{R}}|a_{ijl}(t)|
\cdot\max\limits_{k=1,2}\Big\{|u_{jk}|\Big\}
+\sup\limits_{t\in\mathscr{R}}|J_i(t)|\Big)\Big/\inf\limits_{t\in\mathscr{R}}c_i(t),\,\,\,
\beta_{i1}=\widehat{\ell}_{i1},
_{\left.\begin{ar}{l}\\\end{ar}\right.}\\
\alpha_{i2}=\widetilde{\ell}_{i2}, \,\,\,
\beta_{i2}=\Big(\sum\limits_{l=1}^M\sum\limits_{j=1}^N\sup\limits_{t\in\mathscr{R}}|a_{ijl}(t)|
\cdot\max\limits_{k=1,2}\Big\{|u_{jk}|\Big\}
+\sup\limits_{t\in\mathscr{R}}|J_i(t)|\Big)\Big/\inf\limits_{t\in\mathscr{R}}c_i(t).\nn
\end{ar}\right.
\end{eq}\nid
It is easy for us to check that
$\alpha_{i1}<\beta_{i1}<0<\alpha_{i2}<\beta_{i2}$. The assumption
$(H_1^\mathcal{B})$ implies that $F_i(z)$ is increasing on
$[\ell_{i1}/\sigma_i,\ell_{i2}/\sigma_i]$. Similarly as Theorem 1
and Theorem 2, we have the following two theorems.

\vskip 0.1in\noindent\textsl{Theorem 5}\\
Under the assumptions $(H_1)$ and
$(H_1^\mathcal{B}$)-(H$_2^\mathcal{B}$), each $\mathscr{K}^\Sigma$
is a positively invariant basin with respect to the solution flow
generated by GNNs (1).

\vskip 0.1in\noindent\textsl{Theorem 6}\\Under the basic assumptions
$(H_1)$ and $(H_1^\mathcal{B})$-$(H_2^\mathcal{B})$, for each
$\Sigma$, there exists at least one almost periodic encoded pattern
of GNNs (1) in ${\mathscr{B}}^{\Sigma}$.

\vskip 0.1in\noindent\textsl{Remark 5}\\From $(H_1^\mathcal{B})$ and
piecewise linearity of activation functions in class $\mathcal{B}$ ,
we know that $F_i(z)$ is strictly increasing on
$[\ell_{i1}/\sigma_i,\ell_{i2}/\sigma_i]$ and is strictly decreasing
on $(-\infty,\ell_{i1}/\sigma_i)\cup(\ell_{i2}/\sigma_i,+\infty)$.
By the definition of $\widehat{\ell}_{i1}$,$\widehat{\ell}_{i2}$ and
similar proof of Theorem 1, it is easy for us to know that each
$\mathscr{K}^\Sigma$ is a positively invariant basin with respect to
the solution flow generated by GNNs (1). By Schauder's fixed point
theorem and positive invariancy of each $\mathscr{K}^\Sigma$,
similarly as Theorem 2, we can show that there exists at least one
almost periodic encoded pattern of GNNs (1) in
${\mathscr{B}}^{\Sigma}$.

\vskip 0.1in\indent Since each ${\mathscr{K}}^{\Sigma}$ lies in the
saturated parts to the activation functions of class $\mathcal{B}$,
we get that $\dot{g}_i(z)=0$ for all $z\in[-\infty,
\sigma_i\widehat{\ell}_{i1}]\bigcup[\sigma_i\widetilde{\ell}_{i2},
+\infty]$, that is,
$\mathcal{C}-\mathcal{H}\mathcal{G}\in\mathscr{M}$ always holds. The
exponential stability of almost periodic solutions of GNNs (1)
follows as:

\vskip 0.1in\noindent\textsl{Theorem 7}\\
Under the basic assumptions $(H_1)$ and
$(H_1^\mathcal{B})$-$(H_2^\mathcal{B})$, for each $\Sigma$, there
exists a unique almost periodic encoded pattern of GNNs (1)  which
is exponentially stable in $\mathscr{K}^\Sigma$.

\vskip 0.1in\noindent\textsl{Remark 6}\\ When $c_i$, $a_{ijl}$,
$\kappa_{ijl}$, $J_i:\mathscr{R}\rightarrow\mathscr{R}$ are
$\omega$-periodic functions with $\omega>0$, we also obtain the
existence and exponential stability of $2^N$ periodic solutions of
GNNs (1). For activation functions of class $\mathcal{B}$, we let
$\beta_{i1}:=-\ell_{i1}/\sigma_i$,
$\alpha_{i2}:=\ell_{i2}/\sigma_i$. From Theorem 6 and Theorem 7, we
can prove that each $
\mathscr{K}^\Sigma:=\underbrace{\mathscr{K}_{1\varsigma_1}\times\mathscr{K}_{2\varsigma_2}
\cdots\times\mathscr{K}_{N\varsigma_N}}_N$ is a {\bf larger
attracting basin}.

\vskip 0.15in Now we should consider some special case of GNNs (1)
and compare our results with the existing ones. When $c_i(t)\equiv
c_i$, $M=2$, $\kappa_{ij1}(t)\equiv0$ and $\sigma_j=1$, GNNs (1)
reduces to the following GNNs considered by [20].
\begin{eq}
\frac{du^i(t)}{dt}=\displaystyle-c_iu^i(t)+\sum\limits_{j=1}^Na_{ij1}(t)g_j\big(u^j(t)\big)
+\sum\limits_{j=1}^Na_{ij2}(t)g_j\big(u^j(t-\kappa_{ij2}(t))\big)+J_i(t),
\end{eq}\nid
where $i$$\in$$\mathscr{I}$=$\{1,2,\cdots,N\}$. From Theorem 2 to
Theorem 4, it is easy for us to have the following two corollaries.

\vskip 0.2in\noindent\textsl{Corollary 1}\\
For activation functions of class $\mathcal{A}$, assume  the
following conditions hold:
\begin{ea}
\left\{\begin{ar}{l}\displaystyle (\mathcal{A}_1):\ \
\inf\limits_{\zeta\in\mathscr{R}}\dot{g}_i(\zeta)<
\frac{c_i}{\sigma_i\inf\limits_{t\in\mathscr{R}}\Big[a_{ii1}(t)+a_{ii2}(t)\Big]}<
\sup\limits_{\zeta\in\mathscr{R}}\dot{g}_i(\zeta),
_{\left.\begin{array}{l}\\\end{array}\right.}\\
\displaystyle (\mathcal{A}_2):\ \ (-1)^k\cdot\Big\{F_i(z_{ik})
+J_i(t)\Big\}>\sum\limits_{l=1}^2\sum\limits_{j\neq{i}}\sup\limits_{t\in\mathscr{R}}|a_{ijl}(t)|B_j\
\ (k=1,2),
_{\left.\begin{array}{l}\\\end{array}\right.}\\
(\mathcal{A}_3):\ \ \mathcal{C}-\mathcal{H}\mathcal{G}\in\mathscr{M}
_{\left.\begin{array}{l}\\\end{array}\right.} \mbox{or} \\
\displaystyle (\mathcal{A}_3^*):\ \ pd_ic_i
>\sum\limits_{l=1}^2\sum\limits_{j=1}^N\Big[d_j
\big(\sup\limits_{t\in\mathscr{R}}|a_{jil}(t)|\big)^{pp_{il,m+1}}\big(\dot{g}_i(\xi)\big)^{qq_{il,m+1}}
_{\left.\begin{array}{l}\\\end{array}\right.}\\
\indent\indent\indent+\sum\limits_{k=1}^md_io_k\big(\sup\limits_{t\in\mathscr{R}}|a_{ijl}(t)|\big)^{pp_{jlk}/o_k}
\big(\dot{g}_j(\xi)\big)^{pq_{jlk}/o_k}\Big],
\end{ar}\right.
\end{ea}\nid
where $\mathcal{C}=diag\Big(c_1, c_2,\cdots,c_N\Big)$,
$\dot{g}_j(\xi)=\max\Big\{\dot{g}_j(z)\Big|z=\sigma_j\beta_{j1},\sigma_j\alpha_{j2}\Big\}$;
$z_{ik}$ are defined in Lemma 2; $p>1$, $d_i$, $o_k$, $q_{jlk}$ and
$p_{jlk}$ are positive constants which satisfy with
$\sum\limits_{k=1}^mo_k=p-1$,
$\sum\limits_{k=1}^{m+1}q_{jlk}=\sum\limits_{k=1}^{m+1}p_{jlk}=1$
for each $j\in\mathscr{I}$, $l\in\mathscr{L}$.  Then there exist
only $2^N$ almost periodic encoded patterns of GNNs (18) which are
 exponentially stable.

\vskip 0.1in\noindent\textsl{Corollary 2}\\For activation functions
of class $\mathcal{B}$, assume  the following conditions hold:
\begin{ea}
\left\{\begin{ar}{l}\displaystyle (\mathcal{B}_1):\ \
c_i<\inf\limits_{t\in\mathscr{R}}
\Big[a_{ii1}(t)+a_{ii2}(t)\Big]\sigma_i\dot{\widetilde{g}}_i(\zeta),\,\,\,\zeta\in[\ell_{i1},\ell_{i2}]
_{\left.\begin{array}{l}\\\end{array}\right.}\\
(\mathcal{B}_2):\ \
(-1)^k\cdot\Big\{-c_i\ss\frac{\ell_{ik}}{\sigma_i}+\inf\limits_{t\in\mathscr{R}}\Big[a_{ii1}(t)+a_{ii2}(t)\Big]u_{ik}
+J_i(t)\Big\}
_{\left.\begin{array}{l}\\\end{array}\right.}\\
\indent\indent\indent\ss>\sum\limits_{l=1}^2\sum\limits_{j\neq{i}}
\sup\limits_{t\in\mathscr{R}}|a_{ijl}(t)|\cdot\max\limits_{k=1,2}\Big\{|u_{jk}|\Big\}
\,\,\,(k=1,2).
\end{ar}\right.
\end{ea}\nid
Then there exist only $2^N$ almost periodic encoded patterns of GNNs
(18) which are exponentially stable.

\vskip 0.1in\noindent\textsl{Remark 7}\\For activation functions of
class $\mathcal{A}$, if we only assume that $|g_i(x)-g_i(y)|\leq
\widetilde{L}_i|x-y|$ for all $x,y\in\mathscr{R}$, where
$\widetilde{L}_i$ are positive constants. Let
$\mathcal{G}=diag\Big(\widetilde{L}_1,
\widetilde{L}_2,\cdots,\widetilde{L}_N\Big)$. Under the basic
assumption $(\mathcal{A}_3)$, there exists a unique almost periodic
encoded pattern of GNNs (1) which is globally exponentially stable,
we can refer to [20].

\vskip 0.1in\indent Assume that $c_{i}(t)\equiv c_i$,
$a_{ijl}(t)\equiv a_{ijl}$, $J_i(t)\equiv J_i$, $\sigma_j=1$,
$\kappa_{ijl_{1}}(t)\equiv0$ and $\kappa_{ijl_{2}}(t)\equiv
\kappa_{ijl_{2}}$, where $l_{1}\in\mathscr{L}_{1}$,
$l_{2}\in\mathscr{L}_{2}$ and
$\mathscr{L}_{1}\bigcup\mathscr{L}_{2}=\mathscr{L}$. Then GNNs (1)
reduces to the following autonomous  general neural networks
including [6,8] as our special cases.
\begin{ea}
\frac{du^i(t)}{dt}=\displaystyle-c_iu^i(t)+\sum\limits_{l\in\mathscr{L}_{1}}\sum\limits_{j=1}^Na_{ijl}g_j\big(u^j(t)\big)
+\sum\limits_{l\in\mathscr{L}_{2}}\sum\limits_{j=1}^Na_{ijl}g_j\big(u^j(t-\kappa_{ijl})\big)+J_i,
\end{ea}\nid
where $i$$\in\mathscr{I}$; From Theorem 2 to Theorem 4, it is easy
for us to have the following two corollaries.

\vskip 0.1in\noindent\textsl{Corollary 3}\\For activation functions
of class $\mathcal{A}$, assume  the following conditions hold:
\begin{ea}
\left\{\begin{ar}{l} (\widetilde{\mathcal{A}}_1):\ \
\displaystyle\inf\limits_{\zeta\in\mathscr{R}}\dot{g}_i(\zeta)<
\frac{c_i}{\sum\limits_{l=1}^Ma_{iil}}<
\sup\limits_{\zeta\in\mathscr{R}}\dot{g}_i(\zeta),
_{\left.\begin{array}{l}\\\end{array}\right.}\\
\displaystyle (\widetilde{\mathcal{A}}_2):\ \
(-1)^k\cdot\Big\{F_i(z_{ik})
+J_i\Big\}>\sum\limits_{l=1}^M\sum\limits_{j\neq{i}}|a_{ijl}|B_j\
\ (k=1,2),
_{\left.\begin{array}{l}\\\end{array}\right.}\\
(\widetilde{\mathcal{A}}_3):\ \
\mathcal{C}-\mathcal{H}\mathcal{G}\in\mathscr{M}
_{\left.\begin{array}{l}\\\end{array}\right.} \mbox{or} \\
\displaystyle (\widetilde{\mathcal{A}}_3^*):\ \ pd_ic_i
>\sum\limits_{l=1}^M\sum\limits_{j=1}^N\Big[
d_j|a_{jil}|^{pp_{il,m+1}}\big(\dot{g}_i(\xi)\big)^{qq_{il,m+1}}
+\sum\limits_{k=1}^md_io_k|a_{ijl}|^{pp_{jlk}/o_k}
\big(\dot{g}_j(\xi)\big)^{pq_{jlk}/o_k}\Big],
\end{ar}\right.
\end{ea}\nid
where $\mathcal{C}=diag\Big(c_1, c_2,\cdots,c_N\Big)$,
$\dot{g}_j(\xi)=\max\Big\{\dot{g}_j(z)\Big|z=\sigma_j\beta_{j1},\sigma_j\alpha_{j2}\Big\}$;
$p>1$, $d_i$, $o_k$, $q_{jlk}$ and $p_{jlk}$ are positive constants
which satisfy with $\sum\limits_{k=1}^mo_k=p-1$,
$\sum\limits_{k=1}^{m+1}q_{jlk}=\sum\limits_{k=1}^{m+1}p_{jlk}=1$
for each $j\in\mathscr{I}$, $l\in\mathscr{L}$.  Then there exist
only $2^N$ almost periodic encoded patterns of GNNs (21) which are
 exponentially stable.

\vskip 0.1in\noindent\textsl{Corollary 4}\\For the activation
functions of class $\mathcal{B}$, assume  the following conditions
hold:
\begin{ea}
\left\{\begin{ar}{l}\displaystyle (\widetilde{B}_1):\ \ c_i<
\sum\limits_{l=1}^Ma_{iil}\sigma_i\dot{\widetilde{g}}_i(\zeta),\,\,\,\zeta\in[\ell_{i1},\ell_{i2}]
_{\left.\begin{array}{l}\\\end{array}\right.}\\
\displaystyle (\widetilde{B}_2):\ \
(-1)^k\cdot\Big\{-c_i\frac{\ell_{ik}}{\sigma_i}+\sum\limits_{l=1}^Ma_{iil}u_{ik}
+J_i\Big\}>\sum\limits_{l=1}^M\sum\limits_{j\neq{i}}a_{ijl}\cdot\max\limits_{k=1,2}\Big\{|u_{jk}|\Big\}
_{\left.\begin{array}{l}\\\end{array}\right.}\\
\indent\indent\mbox{where}\,\,\,\, k=1,2.
\end{ar}\right.
\end{ea}\nid
Then there exist only $2^N$ almost periodic encoded patterns of GNNs
(21) which are exponentially stable.

\vskip 0.1in\noindent\textsl{Remark 8}\\In Corollary 3-4, $2^N$
almost periodic encoded patterns of GNNs (21) are indeed equilibria
which are exponential stable. We can replace activation functions of
class $\mathcal{A}$ by
\begin{ea*}
 \left\{\begin{aligned}&g_i\in\mathcal{C}^2,\,\,\,
\eta_i\leq g_i(x)\leq \widetilde{\eta}_i, \,\,\,\dot{g}_i(x)>0,\\
&(x-\vartheta_i)\ddot{g}_i(x)<0, \,\,\mbox{for}
\,\,\,\mbox{all}\,\,\, x\in\mathscr{R},
\end{aligned} \right.
\end{ea*}
where $\eta_i$, $\widetilde{\eta}_i$ and $\vartheta_i$ are constants
with $\eta_i<\widetilde{\eta}_i$, $i\in\mathscr{I}$. There exist
$2^N$ equilibria of GNNs (21) which are exponentially stable. When
$M=2$, $\mathscr{L}_{1}=\{1\}$ and $\mathscr{L}_{2}=\{2\}$, the
related results in [6,8] are the direct results of Corollary 3 and
Corollary 4.  It is obvious that our results are more general than
corresponding results in [6,8].

\vskip 0.1in\noindent\textsl{Remark 9}\\Our approach can also be
adapted to the following general neural networks:
\begin{ea*}
\frac{du^i(t)}{dt}=\displaystyle-c_i(t)u^i(t)+\sum\limits_{l=1}^M
\sum\limits_{j=1}^Na_{ijl}(t)g_j\Big(\sigma_j\int\limits_0^{\kappa_{ijl}}K_{ijl}(s)u^j(t-s)ds\Big)+J_i(t),
\end{ea*}\nid
or
\begin{ea*}
\frac{du^i(t)}{dt}=\displaystyle-c_i(t)u^i(t)+\sum\limits_{l=1}^M
\sum\limits_{j\in
N_l(i)}a_{ijl}(t)g_j\Big(\sigma_ju^j(t-\kappa_{ijl}(t))\Big)+J_i(t),
\end{ea*}\nid
where $K_{ijl}:[0,\kappa_{ijl}]\rightarrow[0,+\infty]$ is assumed be
to continuous and
$0<\int\limits_0^{\kappa_{ijl}}K_{ijl}(s)ds<\infty$,
$N_l(i)=\{i-l,\cdots,i+l\}$, $\kappa_{ijl}\leq+\infty$. The above
general neural networks include [24,38-39,41] as special cases.
Furthermore, our theory generalize stability and existence of
multiple almost periodic (periodic) solutions to above general
neural networks with delays. For more practical applications of
multistability of neural networks, we can refer to [1-3,6-9,48-51].

\vskip 0.3in \normalsize \baselineskip16pt
\centerline{5.~~\textrm{NUMERICAL ILLUSTRATIONS}} \vskip
0.3in\noindent \textsl{Example 1}\\Consider the following neural
networks under almost periodic stimuli.
\begin{eq}
\left\{\begin{ar}{l}\ss\frac{dx_1(t)}{dt}=-(1.2+0.2\cos{2t})x_1(t)
+3g_1(x_1(t))+\sin{\sqrt{7}t}g_2(2x_2(t))
_{\left.\begin{ar}{l}\\\end{ar}\right.}\\
\hskip0.4in+(4+\sin{\sqrt{2}t})g_1(x_1(t))+\cos{\sqrt{3}t}g_2(2x_2(t-9-\sin{t}))+1.1458\cos{\sqrt{5}t},
_{\left.\begin{ar}{l}\\\end{ar}\right.}\\
\ss\frac{dx_2(t)}{dt}=-(3+0.1\sin{3t})x_2(t)
+\cos{\sqrt{3}t}g_1(x_1(t))+4g_2(2x_2(t))
_{\left.\begin{ar}{l}\\\end{ar}\right.}\\
\hskip0.4in+\sin{\sqrt{5}t}g_1(x_1(t-7-3\cos{t}))+(7+\cos{\sqrt{3}t})g_2(2x_2(t))+4.6679\sin{2t},
\end{ar}\right.
\end{eq}\nid
where $g_1(\xi)=g_2(\xi)=\tanh(\xi)$, which belongs to class
$\mathcal{A}$. It is easy for us to get that
\begin{ea*}
F_1(z)=-1.4z+6g(z), \,\,\,\,\,F_2(z)=-3.1z+10g(2z).
\end{ea*}
From some computations, we have $z_{11}=-1.3565$, $z_{12}=1.3565$,
$z_{21}=-0.792$, $z_{22}=0.792$ such that $\dot{F}_{i}(z_{ik})=0$.
Then $F_1(z_{1k})=(-1)^k3.3544$ and $F_2(z_{2k})=(-1)^k6.7370$ where
$k=1,2$. From Lemma 1 and (3)-(4), we can have the following
calculation result:
\begin{ea*}
&&\alpha_{11}=-11.1458,\,\,\,\,\,\beta_{11}=-1.8190,\,\,\,\,\,
\alpha_{12}=1.8190,\,\,\,\,\, \beta_{12}=11.1458,\\
&&\alpha_{21}=-6.4372,\,\,\,\,\,\beta_{21}=-0.9095,\,\,\,\,\,
\alpha_{22}=0.9095,\,\,\,\,\, \beta_{22}=6.4372.
\end{ea*}
It is easy for us to get
\begin{ea*}
\dot{g}_1(\zeta)=\max\Big\{\dot{g}_1(z)\Big|z=\beta_{11},\alpha_{12}\Big\}=0.1,\,\,\,
2\dot{g}_2(\zeta)=\max\Big\{\dot{g}_2(z)\Big|z=2\beta_{21},2\alpha_{22}\Big\}=0.1.
\end{ea*}
Therefore, the parameters satisfy our assumptions in Theorem
3:\\
Assumption $(H_1^{\mathcal{A}})$:
\begin{ea*}
&&0<\sup\limits_{t\in\mathscr{R}}c_1(t)\Big/\inf\limits_{t\in\mathscr{R}}(a_{111}(t)+a_{112}(t))=1.4/6<1,\\
&&0<\sup\limits_{t\in\mathscr{R}}c_2(t)\Big/2\inf\limits_{t\in\mathscr{R}}(a_{221}(t)+a_{222}(t))=3.1/20<1.
\end{ea*}
Assumption $(H_2^{\mathcal{A}})$:
\begin{ea*}
(-1)^k\cdot\{F_1(z_{1k})
+J_1(t)\}\hskip-2mm&=&\hskip-2mm(-1)^k\cdot\{(-1)^k3.3544
+1.1458\cos{\sqrt{5}t}\}\\
\hskip-2mm&>&\hskip-2mm2=\sup\limits_{t\in\mathscr{R}}|a_{121}(t)|B_2+\sup\limits_{t\in\mathscr{R}}|a_{122}(t)|B_2,\\
(-1)^k\cdot\{F_2(z_{2k})
+J_2(t)\}\hskip-2mm&=&\hskip-2mm(-1)^k\cdot\{(-1)^k6.7370
+1.6679\sin{2t}\}\\
\hskip-2mm&>&\hskip-2mm5=\sup\limits_{t\in\mathscr{R}}|a_{211}(t)|B_1+\sup\limits_{t\in\mathscr{R}}|a_{212}(t)|B_1.
\end{ea*}
Assumption $(H_3^{\mathcal{A}})$:
\begin{ea*}
\mathcal{C}-\mathcal{H}\mathcal{G}=\left(\begin{array}{ccc}
1&0\\
0&2.9
\end{array}\right)-\left(\begin{array}{ccc}
8&2\\
2&12
\end{array}\right)\left(\begin{array}{ccc}
0.1&0\\
0&0.1
\end{array}\right)=\left(\begin{array}{ccc}
0.2&-0.2\\
-0.2&1.7
\end{array}\right)\in\mathscr{M}.
\end{ea*}
\begin{figure}
\begin{center}
\resizebox{!}{90mm}{\includegraphics*[56,193][560,600]{newm1.eps}}
\end{center}
\caption{Convergence dynamics of four almost periodic encoded
patterns of (24).}
\end{figure}Then there exist four almost periodic encoded patterns of (24) in
$\mathscr{B}^\Sigma$ and their attracting basins are
$\mathscr{K}^\Sigma$. We can compute the following regions:\vskip
0.1in Region I: $(-\infty,
-1.77]\times[0.9,+\infty)$,\,\,\,\,\,\,~~~ Region II: $[1.77,
+\infty]\times[0.9,+\infty)$,\vskip 0.1in

Region III: $(-\infty, -1.77]\times(-\infty, -0.9]$,\,\,\, Region
IV: $[1.77, +\infty]\times(-\infty, -0.9]$\vskip 0.1in\noindent
which contain larger attracting basins mentioned in Remark 3. From
Figure 1, we can see that four almost periodic encoded patterns of
(24) lie in invariant regions $\mathscr{B}^{(1,2)}$,
$\mathscr{B}^{(2,2)}$, $\mathscr{B}^{(1,1)}$ and
$\mathscr{B}^{(2,1)}$  which borderlines are plotted in blue. Their
larger attracting basins are denoted by Region I to Region IV which
borderlines are plotted in black.

\vskip 0.1in\noindent\textsl{Example 2}\\Consider the following
neural networks under periodic stimuli.
\begin{eq}
\left\{\begin{ar}{l}\ss\frac{dx_1(t)}{dt}=-1.4x_1(t)
+2g_1(x_1(t))+0.1g_2(2x_2(t))
_{\left.\begin{ar}{l}\\\end{ar}\right.}\\
\hskip0.4in+4g_1(x_1(t-10))+0.1g_2(2x_2(t-10))+3.1456\sin{t},
_{\left.\begin{ar}{l}\\\end{ar}\right.}\\
\ss\frac{dx_2(t)}{dt}=-3.1x_2(t) +0.1g_1(x_1(t))+3g_2(2x_2(t))
_{\left.\begin{ar}{l}\\\end{ar}\right.}\\
\hskip0.4in+0.1g_1(x_1(t-10))+7g_2(2x_2(t-10))+6.1705\cos{t},
\end{ar}\right.
\end{eq}\nid
where $g_1(\xi)=g_2(\xi)=\tanh(\xi)$, which belongs to class
$\mathcal{A}$. Similarly as Example 1, we can check that
$(H_1^{\mathcal{A}})$ and $(H_2^{\mathcal{A}})$ hold. From Lemma 1
and (3)-(4), we can have the following calculation result:
\begin{ea*}
&&\alpha_{11}=-6.6754,\,\,\,\,\,\beta_{11}=-1.443,\,\,\,\,\,
\alpha_{12}=1.443,\,\,\,\,\, \beta_{12}=6.6754,\\
&&\alpha_{21}=-5.2808,\,\,\,\,\,\beta_{21}=-1.0891,\,\,\,\,\,
\alpha_{22}=1.0891,\,\,\,\,\, \beta_{22}=5.2808,\\
&&\dot{g}_1(\zeta)=\max\Big\{\dot{g}_1(z)\Big|z=\beta_{11},\alpha_{12}\Big\}=0.2,\,\,\,
2\dot{g}_2(\zeta)=\max\Big\{2\dot{g}_2(z)\Big|z=2\beta_{21},2\alpha_{22}\Big\}=0.1.
\end{ea*}
Let $p=4$, $m=1$, $o_1=3$, $q_{ilk}=p_{ilk}=1/2$ ($i,l,k=1,2$). From
some calculations, we can check that $(H_4^{\mathcal{A}})$ holds.
Therefore, by Theorem 4, there exist four periodic encoded patterns
of (25) in $\mathscr{B}^\Sigma$. Their attracting basins follow
as:\vskip 0.1in Region I: $(-\infty,
-1.443]\times[1.0891,+\infty)$,\,\,\,~~~~~ Region II: $[1.443,
+\infty]\times[1.0891,+\infty)$,\vskip 0.1in

Region III: $(-\infty, -1.443]\times(-\infty, -1.0891]$,\,\,\,
Region IV: $[1.443, +\infty]\times(-\infty, -1.0891]$.\vskip
0.1in\noindent We can refer to their convergence dynamics plotted in
Figure 2.
\begin{figure}
\begin{center}
\resizebox{!}{90mm}{\includegraphics*[56,193][560,600]{newp1.eps}}
\end{center}
\caption{Convergence dynamics of four periodic encoded patterns of
(25).}
\end{figure}
\begin{figure}
\begin{center}
\resizebox{!}{90mm}{\includegraphics*[56,193][560,600]{newpm1.eps}}
\end{center}
\caption{Convergence dynamics of four almost periodic encoded
patterns of (26).}
\end{figure}
\vskip 0.1in\noindent\textsl{Example 3}\\Consider the following
neural networks under almost periodic stimuli.\\
\begin{eq}
\left\{\begin{ar}{l}\ss\frac{dx_1(t)}{dt}=-(1.5+0.5\cos{\sqrt{7}t})x_1(t)
+(4.5+0.5\sin{\sqrt{2}t})g_1(0.5x_1(t))+0.1\sin{t}g_2(0.25x_2(t))
_{\left.\begin{ar}{l}\\\end{ar}\right.}\\
\hskip0.4in+(2.5+0.5\sin{\sqrt{2}t})g_1(0.5x_1(t-10))+0.1\cos{t}g_2(0.25x_2(t-10))+0.2\sin{3t},
_{\left.\begin{ar}{l}\\\end{ar}\right.}\\
\ss\frac{dx_2(t)}{dt}=-(0.75+0.25\sin{t})x_2(t)
+0.2\cos{\sqrt{5}t}g_1(0.5x_1(t))+(8+\cos{\sqrt{3}t})g_2(0.25x_2(t))
_{\left.\begin{ar}{l}\\\end{ar}\right.}\\
\hskip0.4in+0.2\sin{\sqrt{7}t}g_1(0.5x_1(t-10))+(2+\cos{\sqrt{3}t})g_2(0.25x_2(t-10))+0.05\cos{2t},
\end{ar}\right.
\end{eq}\nid
where $g_1(\xi)=g_2(\xi)=\frac{1}{2}(|\xi+1|+|\xi-1|)$, which
belongs to class $\mathcal{B}$. Similarly as Example 1, we can check
that $(H_1^{\mathcal{B}})$ and $(H_2^{\mathcal{B}})$ hold. From
Lemma 1 and (3)-(4), we can have the following calculation result:
\begin{ea*}
&&\alpha_{11}=-8.4,\,\,\,\,\,\beta_{11}=-2.8,\,\,\,\,\,
\alpha_{12}=2.8,\,\,\,\,\, \beta_{12}=8.4,\\
&&\alpha_{21}=-24.9,\,\,\,\,\,\beta_{21}=-7.55,\,\,\,\,\,
\alpha_{22}=7.55,\,\,\,\,\, \beta_{22}=24.9.
\end{ea*}
By Theorem 7, there exist four stable almost periodic encoded
patterns of (26) in $\mathscr{B}^\Sigma$. It is obvious that their
larger attracting basins follow as:\vskip 0.1in Region I: $(-\infty,
-2]\times[4,+\infty)$,\,\,\, ~~~~~Region II: $[2,
+\infty]\times[4,+\infty)$,\vskip 0.1in

Region III: $(-\infty, -2]\times(-\infty, -4]$,\,\,\, Region IV:
$[2, +\infty]\times(-\infty, -4]$.\vskip 0.1in\noindent Their
convergence dynamics are illustrated in Figure 3.

\vskip 0.1in\noindent\textsl{Example 4}\\Consider the following
neural networks under periodic stimuli.
\begin{eq}
\left\{\begin{ar}{l}\ss\frac{dx_1(t)}{dt}=-(1.5+0.5\cos{t})x_1(t)
+10g_1(0.5x_1(t))+0.1g_2(0.25x_2(t))
_{\left.\begin{ar}{l}\\\end{ar}\right.}\\
\hskip0.4in+3g_1(0.5x_1(t-10))+0.1g_2(0.25x_2(t-10))+0.2\sin{4t},
_{\left.\begin{ar}{l}\\\end{ar}\right.}\\
\ss\frac{dx_2(t)}{dt}=-(0.8+0.2\sin{t})x_2(t)
+0.2g_1(0.5x_1(t))+9g_2(0.25x_2(t))
_{\left.\begin{ar}{l}\\\end{ar}\right.}\\
\hskip0.4in+0.2g_1(0.5x_1(t-10))+3g_2(0.25x_2(t-10))+0.05\cos{8t},
\end{ar}\right.
\end{eq}\nid
where $g_1(\xi)=g_2(\xi)=\frac{1}{2}(|\xi+1|+|\xi-1|)$, which
belongs to class $\mathcal{B}$. Similarly as Example 1, we can check
that $(H_1^{\mathcal{B}})$ and $(H_2^{\mathcal{B}})$ hold. By some
calculations, we can have the following result:
\begin{ea*}
&&\alpha_{11}=-13.4,\,\,\,\,\,\beta_{11}=-6.3,\,\,\,\,\,
\alpha_{12}=6.3,\,\,\,\,\, \beta_{12}=13.4,\\
&&\alpha_{21}=-20.75,\,\,\,\,\,\beta_{21}=-11.55,\,\,\,\,\,
\alpha_{22}=11.55,\,\,\,\,\, \beta_{22}=20.75.
\end{ea*}
By Theorem 7, there exist four stable periodic encoded patterns of
(27) in $\mathscr{B}^\Sigma$. Their larger attracting basins follow
as:\vskip 0.1in Region I: $(-\infty, -2]\times[4,+\infty)$,\,\,\,
Region II: $[2, +\infty]\times[4,+\infty)$,\vskip 0.1in

Region III: $(-\infty, -2]\times(-\infty, -4]$,\,\,\, Region IV:
$[2, +\infty]\times(-\infty, -4]$.\vskip 0.1in\noindent For their
corresponding convergence dynamics, we can refer to Figure 4.
\begin{figure}
\begin{center}
\resizebox{!}{90mm}{\includegraphics*[56,193][560,600]{newpp1.eps}}
\end{center}
\caption{Convergence dynamics of four periodic encoded patterns of
(27).}
\end{figure}
\vskip 0.3in \normalsize \baselineskip16pt
\centerline{6.~~\textrm{CONCLUDING REMARKS}} \vskip 0.15in\noindent
In this paper, we investigate multi-almost periodicity of general
neural networks under almost periodic stimuli. Invariant regions and
attracting basins are established to investigate existence and
exponential stability of $2^N$ almost periodic encoded patterns. Our
results extend and generalize the related results reported in the
literature [6,8,10,20-23].

\vskip 0.3in \normalsize \baselineskip16pt
\centerline{\textrm{REFERENCES}} \vskip 0.1in\footnotesize
\baselineskip12pt

\noindent 1. L. O. Chua, L. Yang. Cellular neural networks: Theory.
\textit{IEEE Trans. Circuits Syst.} 1988; \textbf{35}:1257-1272.\\
2. L. O. Chua, L. Yang. Cellular neural networks: Application.
\textit{IEEE Trans. Circuits Syst.}  1988; \textbf{35}:1273-1290.\\
3. J. Foss, A. Longtin, B. Mensour, J.Milton. Multistability and
delayed recurrent loops. \textit{Phys. Rev. Lett.} 1996;
\textbf{76}:708-711.\\
4. J. Hopfield. Neurons with graded response have collective
computational properties like those of two state neurons.
\textit{Proc. Natl. Acad. Sci. USA} 1984; \textbf{81}:3088-3092.\\
5. F. Forti. On global asymptotic stability of a class of nonlinear
systems arising in neural network theory. \textit{J. Diff. Eqs.}
1994; \textbf{113}:246-264.\\
6. C. Y. Cheng, K. H. Lin, C. W. Shih. Multistability in recurrent
neural networks. \textit{SIAM J. Appl. Math.} 2006;
\textbf{66}(4):1301-1320.\\
7. Z. Zeng, D. S. Huang, Z. Wang. Memory pattern analysis of
cellular neural networks. \textit{Phys. Lett. A} 2005;
\textbf{342}:114-128.\\
8. C. Y. Cheng, K. H. Lin, C. W. Shih. Multistability and
convergence in delayed  neural networks. \textit{Physica D} 2007;
\textbf{225}:61-74.\\
9. L. P. Shayer, S. A. Campbell. stability, bifurcation and
multistability in a system of two coupled neurons with multiple time
delays. \textit{SIAM J. Appl. Math.} 2000; \textbf{61}:673-700.\\
10. Mohamad S. Convergence dynamics of delayed Hopfield-type neural
networks under almost periodic stimuli. \textit{Acta Appl. Math.}
2003; \textbf{76}:117-135.\\
11. K. Gopalsamy. Stability and oscillations in delay differential
equations of population dynamics. Kluwer Academic Publishers, The
Netherlands, 1992.\\
12. A. M. Fink. Almost periodic differential equations. Springer,
New York, 1974.\\
13. C. Y. He. Almost periodic differential equations. Beijing,
China: Higher Education, 1992.\\
14. Lu W. L., Chen T. P. Global stability of almost periodic
solution for delayed neural networks. \textit{Chin. Sci. A} 2005;
\textbf{35}(7):774-784.\\
15. Huang ZK, Wang XH, Gao F. The existence and global attractivity
of almost periodic sequence solution of discrete-time neural
networks. \textit{Phys Lett A} 2006;\textbf{350}:182-191.\\
16. Huang ZK, Xia YH, Wang XH. The existence and exponential
attractivity of $\kappa$-almost periodic sequence solution of
discrete time neural
networks. \textit{Nonlinear Dynamics} 2007;\textbf{50}:13-26.\\
17. Takhashi N. A new sufficient condition for complete stability of
cellular neural networks with delay. \textit{IEEE Trans. Circuits
Syst. I} 2000; \textbf{47}(6):793-799.\\
18. Chunhua Feng, R. Plamondon. On the stability analysis of delayed
neural networks systems. \textit{Neural Networks} 2001;
\textbf{14}:1181-1188.\\
19. Chunhua Feng, Peiguang Wang. Almost periodic solutions of forced
lienard-type equations with time delays. \textit{J. Comput. Appl.
Math.} 2003; \textbf{161}:67-74.\\
20. B. W. Liu, L. H. Huang. Existence and exponential stability of
almost periodic solutions for cellular neural networks with
time-varying delays. \textit{Phys. Lett. A}  2005;
\textbf{341}:135-144.\\
21. H. Y. Zhao. Existence and global exponential convergence of
almost periodic solutions for cellular neural networks with variable
delay. \textit{J. Eng. Math.}   2005; \textbf{22}(2):295-300.\\
22. H. Y. Zhao, G. L. Wang. Existence and global attractivity of
almost periodic solutions for Hopfield neural networks with variable
delay. \textit{Math. Acta Scientia}  2004; \textbf{24}(6):723-729.\\
23. A. P. Chen, L. H. Huang. Existence and attractivity of almost
periodic solutions of Hopfield neural networks. \textit{Math. Acta
Scientia}  2001; \textbf{21}(4):505-511.\\
24. H. Huang, Daniel W.C. Ho, J. D. Cao. Analysis of global
exponential stability and periodic solutions of neural networks with
time-varying delays. \textit{Neural Networks}  2005;
\textbf{18}:161-170.\\
25. J. D. Cao, Q. Li. On the exponential stability and periodic
solutions of delayed cellular neural networks. \textit{J. Math.
Anal. Appl.}  2000; \textbf{252}:50-64.\\
26. H. Fang, J. B. Li. Global exponential stability and periodic
solutions of cellular neural networks with delay. \textit{Phys. Rev.
E}  2000; \textbf{61}(4):4212-4217.\\
27. J. D. Cao. Periodic solutions and exponential stability in
delayed cellular neural networks. \textit{Phys. Rev. E}
 1999; \textbf{60}(3):3244-3248.\\
28. Cao JD, Wang L. Exponential stability and periodic oscillatory
solution in BAM networks with delays.  \textit{IEEE Trans Neural
Networks}  2002; \textbf{13}:457-463.\\
29. Song QK, Wang ZD. An analysis on existence and global
exponential stability of periodic solutions for BAM neural networks
with time-varying delays. \textit{Nonlinear Analysis: Real World
Applications} 2006; doi:10.1016/j.nonrwa.2006.07.002.\\
30. Chen AP,Huang LH, Liu ZG, Cao JD. Periodic bidirectional
associative memory neural networks with distributed delays.
\textit{J Math Anal Appl}  2006; \textbf{317}:80-102.\\
31. H. J. Jiang, Z. D. Teng. Dynamics of neural networks with
variable coefficients and time-varying delays. \textit{Neural
Networks}  2006; \textbf{19}:676-683.\\
32. Berman A, Plemmons RJ.  Nonnegative Matrices in the Mathematical
Science. Academic Press New York 1929.\\
33. R. A. Horn, C. R. Johnson. Topics in Matrix Analysis. Cambridge
Univ Press Cambridge UK 1991.\\
34. K. N. Lu, D. Y. Xu, Z. C. Yang. Global attraction and stability
for Cohen-Grossberg neural networks with delays. \textit{Neural
networks} doi:10.1016/j.neunet.2006.07.006.\\
35. Gopalsamy K., X. Z. He. Delay-independent stability in
bi-directional associative memory networks. \textit{IEEE Trans.
Neural Networks}  1994; \textbf{5}:998-1002.\\
36. Liao XF, Yu JB. Qualitative analysis of bidrectional associative
memory with time delays.  \textit{Int. J. Circuit Theory Appl.}
 1998; \textbf{26}:219-229.\\
37. L. H. Huang, C. X. Huang, B. W. Liu. Dynamics of a class of
cellular neural networks with time-varying delays. \textit{Phys.
Lett. A}  2005; \textbf{345}:330-344.\\
38. T. Roska, L. O. Chua. Cellular neural networks with non-linear
and delay-type template elements and non-uniform grids. \textit{Int.
J. Circuit Theory Appl.}  1992; \textbf{20}:469-481.\\
39. Y. G. Liu, Z. S. You, L. P. Cao. On the almost periodic solution
of cellular neural networks with distributed delays. \textit{IEEE
Trans. Neural Networks}  2007; \textbf{18}:295-300.\\
40. M. Gilli. Stability of cellular neural networks and delayed
cellular neural networks with nonpositive templates and nonmonotonic
output functions. \textit{IEEE Trans. Circ. sys. I} 1994; \textbf{41}:518-528.\\
41. S. J. Guo, L. H. Huang. Exponential stability and periodic
solutions of neural networks with continuously distributed delays.
\textit{Physical Review E} 2003; \textbf{67}:011902.\\
42. M. Gilli. Analysis of periodic oscillations in
finite-dimensional CNNs through a spatio-temporal harmonic balance
technique. \textit{Int. J. Circ. Theor. Appl.} 1998;
\textbf{25}(4):279-288.\\
43. V. Lanza, F. Corinto, M. Gilli, Pier P. Civalleri. Analysis of
nonlinear oscillatory network dynamics via time-varying amplitude
and phase variables.
 \textit{Int. J. Circ. Theor. Appl.} 2007; \textbf{35}(5-6):623-644.\\
44. P. P. Civalleri, M. Gilli, M. Bonnin. Equivalent circuits for
small signal performance of spin $1/2$ particles
. \textit{Int. J. Circ. Theor. Appl.} 2006; \textbf{34}(2):165-182.\\
45. I. Petr$\acute{a}$s, M. Gilli. Complex dynamics in
one-dimensional CNNs.
 \textit{Int. J. Circ. Theor. Appl.} 2006; \textbf{34}(1):3-20.\\
46. M. Forti. M-matrices and global convergence of discontinuous
neural networks.
 \textit{Int. J. Circ. Theor. Appl.} 2007; \textbf{35}(2):105-130.\\
47. Huang ZK, Xia YH, Wang XH. Exponential stability of impulsive
Cohen-Grossberg networks with distributed delays.
 \textit{Int. J. Circ. Theor. Appl.}, in press,
 doi:10.1002/cta.424.\\
48. M. Brucoli, L. Carnimeo, G. Grassi. Discrete-time cellular
neural networks for associative memories with learning and
forgetting capabilities.
 \textit{IEEE Trans. Circ. Sys. I} 1995; \textbf{42}(7):396-399.\\
49. J. Juang, S. S. Lin. Cellular neural networks: Mosaic pattern
and spatial chaos. \textit{SIAM J. Appl. Math.} 2000; \textbf{60}:
891-915.\\
50. C. W. Shih, Pattern formation and spatial chaos for cellular
neural networks with asymmetric templates. \textit{Internat. J.
Bifur. Chaos Appl. Sci. Engrg.} 1998; \textbf{8}:1907-1936.\\
51. C. W. Shih. Influence of boundary conditions on pattern
formation and spatial chaos in lattice systems. \textit{SIAM J.
Appl. Math.} 2000; \textbf{61}:335-368.
\end{document}